\documentclass[11pt]{amsart}
\usepackage{enumerate, amsfonts, amsmath, amssymb, amsthm, mathtools, thmtools, wasysym, graphics, graphicx, xcolor, frcursive,xparse,comment,ytableau, import,dynkin-diagrams}
\usepackage[margin=1in]{geometry}
\usepackage[vcentermath]{youngtab}
\usepackage{genyoungtabtikz}
\usepackage{ mathrsfs }
\usepackage[all]{xy}
\usepackage{url, hyperref, hypcap}
\hypersetup{colorlinks=true, citecolor=darkblue, linkcolor=darkblue}

\usepackage{tikz}
\usetikzlibrary{calc,through,backgrounds,shapes,matrix}

\definecolor{darkblue}{rgb}{0.0,0,0.7} 
\newcommand{\darkblue}{\color{darkblue}} 
\definecolor{darkred}{rgb}{0.7,0,0} 
\definecolor{lightgrey}{rgb}{0.7,0.7,0.7} 
\makeatletter

\usepackage{graphicx}
\def\chk#1{#1^{\normalfont\smash{\scalebox{.7}[1.4]{\rotatebox{90}{\guilsinglleft}}}}}
\def\Z{\mathbb{Z}}
\def\R{\mathbb{R}}

\def\simp{\Delta}
\def\wa{\widetilde{W}}
\def\ac{\mathcal{A}}

\def\x{\mathbf{x}}

\def\h{\mathrm{ht}}
\def\cw{\chk{\omega}}

\def\rhoh{\frac{\chk{\rho}}{h}}

\def\rat{b}
\def\Q{\mathcal{Q}}
\def\bt{t_b}
\def\br{r_b}

\newcommand{\cwl}{\chk{\Lambda}}
\newcommand{\wex}{\wa_{\mathrm{ex}}}
\newcommand{\cycl}{\Omega}

\newcommand{\Expt}[2]{\operatorname*{\mathbb{E}}\limits_{#1}(#2)}

\renewcommand{\max}[2]{\operatorname*{max}\limits_{#1}(#2)}
\newcommand{\core}{{\sf core}}

\newcommand{\size}{{\sf size}}
\newcommand{\sizeb}{\size^{(b)}}

\newcommand{\Sommers}{\mathcal{S}}

\newcommand{\lcm}{{\sf lcm}}

\newcommand{\inv}{{\sf inv}}

\renewcommand{\mod}{\operatorname{mod}}

\newcommand{\waf}{\widetilde{w}}
\newcommand{\wf}{\waf_{\rat}}

\newcommand{\amax}{\tilde{\alpha}}
\newtheorem{theorem}{Theorem}[section]
\newtheorem{proposition}[theorem]{Proposition}

\theoremstyle{definition}
\newtheorem{definition}[theorem]{Definition}
\newtheorem{example}[theorem]{Example}
\newtheorem{conjecture}[theorem]{Conjecture}

\newtheorem{remark}[theorem]{Remark}
\newtheorem{convention}[theorem]{Convention}

\usepackage[nameinlink]{cleveref}
\usepackage[all]{xy}
\usepackage[T1]{fontenc}
\crefformat{footnote}{#2\footnotemark[#1]#3}
\crefformat{conjecture}{Conjecture~#2#1#3}
\crefformat{convention}{Convention~#2#1#3}

\newcommand{\defn}[1]{\emph{\darkblue #1}}
\usetikzlibrary{math}
\usepackage{xifthen}
\usepackage{xstring}
\newcommand\IfStringInList[2]{\IfSubStr{,#2,}{,#1,}}

\title{Strange Expectations in Affine Weyl Groups}

\author[E.N.~Stucky]{Eric Nathan Stucky}
\address[E.N.~Stucky]{University of Minnesota-Twin Cities}
\address[E.N.~Stucky]{Charles University, Faculty of Mathematics and Physics, Department of Algebra, Sokolovska 83, 18600 Praha 8, Czech Republic}
\email{stuck127@umn.edu}

\author[M.~Thiel]{Marko Thiel}
\address[M.~Thiel]{Unaffiliated}
\email{thiel.marko@gmail.com}

\author[N.~Williams]{Nathan Williams}
\address[N.~Williams]{University of Texas at Dallas}
\email{nathan.f.williams@gmail.com}

\date{\today}

\begin{document}

\begin{abstract}
Our main result is a generalization, to all affine Weyl groups, of P. Johnson's proof of D. Armstrong's conjecture for the expected number of boxes in a simultaneous core. This extends earlier results by the second and third authors in simply-laced type. We do this by modifying and refining the appropriate notion of the "size" of a simultaneous core. In addition, we provide combinatorial core-like models for the coroot lattices in classical type and type $G_2$.
\end{abstract}

\maketitle

\section{Introduction}
\label{sec:introduction}

\subsection{Motivation}

Macdonald's celebrated affine denominator formula \[\prod_{\alpha \in \widetilde{\Phi}^+} (1-e^{-\alpha})^{\mathrm{mult}(\alpha)} = \sum_{w \in \widetilde{W}} (-1)^{\ell(w)} e^{w(\rho)-\rho}\]
specializes to many famous identities, including 
 Euler's pentagonal number theorem,
 Jacobi's triple product identity,
 and Dyson's identity for Ramanujan's $\tau$-function~\cite{macdonald1971affine,dyson1972missed}. One such specialization---for simply-laced types---is the equality
\begin{equation}\prod_{i=1}^\infty c(x^i) = \left(\prod_{i=1}^\infty \frac{1}{1-x^{hi}} \right)^n \sum_{q \in \mathcal{\Q}} x^{\left\langle \frac{h}{2}q - \rho,q\right\rangle},\label{eq:mac}\end{equation} where
$h$ is the \defn{Coxeter number} and
$c(x)$ is the characteristic polynomial of a \defn{Coxeter element}.  There is a version for all types, which Macdonald refers to but omits in~\cite{macdonald1971affine}\footnote{At the end of~\cite[Section 8]{macdonald1971affine}, Macdonald writes ``When $R$ contains roots of different lengths, the formula corresponding to [\Cref{eq:mac}] is more complicated, and we shall not reproduce it here.''}:
\begin{equation}\sum_{q \in \mathcal{\Q}} x^{\left\langle \frac{h}{2}q - \rho,q\right\rangle} = \prod\limits_{i=1}^\infty\left[ (1-x^i)^{n_s}(1-x^{ri})^{n_\ell} \left(\prod\limits_{\alpha \in \Phi_s} (1-x^i \omega^{\mathrm{ht}(\alpha)})\right)\left(\prod\limits_{\alpha \in \Phi_\ell} (1-x^{ri} \omega^{\mathrm{ht}(\alpha)})\right)\right],\label{eq:thiel}\end{equation}
where
$n_s$ and $n_\ell$ count the number of short and long roots,
$\omega$ is a primitive $h$th root of unity,
$r$ is the ratio of the length of a long to short root,
$\Phi_s$ and $\Phi_\ell$ are the sets of short and long roots, and
$\text{ht}(\alpha)$ is the height of the root $\alpha$.

\subsection{Partitions}
Recall that an \defn{integer partition} is a sequence of non-increasing positive integers $\lambda=(\lambda_1 \geq \lambda_2 \geq \cdots \geq \lambda_k)$.  The \defn{Ferrers diagram} of an integer partition $\lambda$ (under the English convention) is a top-left justified subset of $\mathbb{N}\times\mathbb{N}$ with $\lambda_i$ boxes in the $i$-th row (counting from the top).  A \defn{hook} of a given box in a Ferrers diagram is the collection of boxes to the right and below the given box.  An example is given in~\Cref{fig:partition}.
\begin{figure}[htbp]
\[\begin{ytableau}
8 & 5 & 4 & 2 & 1 \\
5 & 2 & 1 \\
2\\
1
\end{ytableau}\]
\caption{The partition $\lambda=(5,3,1,1)$, with boxes labeled by the lengths of their hooks.  It is a $3$-core (since it has no hooks of size $3$) but not a $4$-core.}
\label{fig:partition}
\end{figure}

An \defn{$a$-core} is an integer partition with no hook of length $a$.  For example, in type $A$,~\Cref{eq:mac} can be interpreted as the beautiful combinatorial formula 
\begin{equation}\label{eq:typea}\prod_{i=1}^\infty \frac{1}{1-x^i} = \left(\prod_{i=1}^\infty \frac{1}{1-x^{ai}}\right)^a\sum_{q \in \core(a)} x^{\size(q)}.\end{equation}

An \defn{$(a,b)$-core} is a partition that is simultaneously an $a$-core and a $b$-core.  For $a$ and $b$ relatively prime, it turns out that there are only finitely many $(a,b)$-cores: \[\big|\core(a,b)\big| = \frac{1}{a+b}\binom{a+b}{b}.\]
For $\lambda$ a partition, write $\lambda^\intercal$ for its conjugate and $\size(\lambda)$ for the number of its boxes.  The starting point for a number of recent investigations has been Armstrong's conjecture on the average number of boxes in an $(a,b)$-core, and in a self-conjugate $(a,b)$-core~\cite{armstrong2011conjecture,armstrong2014results}, which can be thought of as a sort of finite version of~\Cref{eq:typea}.

\begin{theorem}[{\cite{johnson2015lattice}}]\label{thm:arm}  For $\gcd(a,b)=1$,
\begin{equation*}
\Expt{\lambda \in \core(a,b)}{\size(\lambda)} = \frac{(a-1)(b-1)(a+b+1)}{24} = \Expt{\substack{\lambda \in \core(a,b) \\ \lambda = \lambda^\intercal}}{\size(\lambda)}.
\end{equation*}
\end{theorem}
Both equalities in~\Cref{thm:arm} were proven by Johnson using weighted Ehrhart theory~\cite{johnson2015lattice}; the second equality was first proven by Chen, Huang, and Wang~\cite{chen2014average}. 

In~\cite{thiel2017strange}, we generalized Armstrong's conjecture and Johnson's proof of the first equality to all \emph{simply-laced} affine Weyl groups, thereby giving a sort of finite analogue of~\Cref{eq:mac}. In the present paper, we find and---in the words of Macdonald---\emph{reproduce} the generalization to all affine Weyl groups, giving a finite version of~\Cref{eq:thiel}.

\subsection{Combinatorial Models of Coroot Lattices}

The set of $a$-cores under the action of the affine symmetric group $\widetilde{\mathfrak{S}}_a$ is a well-studied combinatorial model for the coroot lattice $\chk{\Q}_{a}$ of type $A_{a-1}$. Indeed, for all affine Weyl groups $\wa=\wa(X_n) := W \ltimes \chk{\Q}_{X_n} $, there is a well-known $\wa$-equivariant map from the group to the coroot lattice $\wa \to \chk{\Q}_{X_n}$ given by $\waf \mapsto \waf(0)$, which restricts to a $\wa$-equivariant bijection on the cosets $\wa/W$. Thus, combinatorial models for $\chk{\Q}_{X_n}$ also give models for $\wa/W$, representatives usually taken to be dominant affine elements.  In type $A_{a-1}$, these correspondences give $\widetilde{\mathfrak{S}}_a$-equivariant bijections
\begin{equation}
\begin{alignedat}{2}
\core(a) &\leftrightarrow \chk{\Q}_{A_{a-1}} &&\leftrightarrow \widetilde{\mathfrak{S}}_a / \mathfrak{S}_{a} \\
 \lambda &\leftrightarrow q_\lambda &&\leftrightarrow \waf_\lambda.
\end{alignedat}
\label{eq:core_to_coroot}
\end{equation}
We describe the first of these bijections ($\lambda\leftrightarrow \chk{\Q}_{A_{a-1}}$) in detail in \Cref{sec:abaci}.

\begin{remark}
Here and throughout, all actions of $\wa$ on various sets are left-actions. In particular, this means that given a reduced expression $a_1\cdots a_k$, the corresponding simple transpositions act in \emph{decreasing} order of the indices (that is, we ``read'' right to left).
\end{remark}

\begin{figure}[htbp]
\begin{center}


\resizebox{.32\textwidth}{!}{
\begin{tikzpicture}[scale =1,bull/.style={circle,draw=black,fill=black,thick,scale=.6}]
\tikzmath{\xpoint=1.732050808;\ypoint=1;\xbdry=1.154700538;\ybdry=1;\hght=8;}
\clip (.8660254038,.5) circle (3.1);
\newcounter{co};
\newcounter{li};
\newcounter{xcurr};
\newcounter{ycurr};

\setcounter{li}{-2}
\whiledo{\value{li}<3}
{
    \setcounter{co}{-3}
    \whiledo{\value{co}<4}
    {
        \setcounter{xcurr}{2*\value{li}};
        \setcounter{ycurr}{2*\value{co}};
        \node at (\value{xcurr}*\xpoint,\value{ycurr}*\ypoint)[bull]{};
        \draw[-,black,very thin](\value{xcurr}*\xpoint+\xbdry,\value{ycurr}*\ypoint)--(\value{xcurr}*\xpoint+.5*\xbdry,\value{ycurr}*\ypoint+\ybdry)--(\value{xcurr}*\xpoint-.5*\xbdry,\value{ycurr}*\ypoint+\ybdry)--(\value{xcurr}*\xpoint-\xbdry,\value{ycurr}*\ypoint)--(\value{xcurr}*\xpoint-.5*\xbdry,\value{ycurr}*\ypoint-\ybdry)--(\value{xcurr}*\xpoint+.5*\xbdry,\value{ycurr}*\ypoint-\ybdry)--cycle;
        \stepcounter{co};
    }
    \stepcounter{li};
}

\setcounter{li}{-2};
\whiledo{\value{li}<2}
{
    \setcounter{co}{-4};
    \whiledo{\value{co}<4}
    {
        \setcounter{xcurr}{2*\value{li}+1};
        \setcounter{ycurr}{2*\value{co}+1};
        \node at (\value{xcurr}*\xpoint,\value{ycurr}*\ypoint)[bull]{};
        \draw[-,black,very thin](\value{xcurr}*\xpoint+\xbdry,\value{ycurr}*\ypoint)--(\value{xcurr}*\xpoint+.5*\xbdry,\value{ycurr}*\ypoint+\ybdry)--(\value{xcurr}*\xpoint-.5*\xbdry,\value{ycurr}*\ypoint+\ybdry)--(\value{xcurr}*\xpoint-\xbdry,\value{ycurr}*\ypoint)--(\value{xcurr}*\xpoint-.5*\xbdry,\value{ycurr}*\ypoint-\ybdry)--(\value{xcurr}*\xpoint+.5*\xbdry,\value{ycurr}*\ypoint-\ybdry)--cycle;
        \stepcounter{co};
    }
    \stepcounter{li};
}
\setcounter{co}{0}
\whiledo{\value{co}<8}
{
    \draw[-,black,thick] (\value{co}*\xbdry,0)--(\value{co}*\xbdry+.5*\hght*\xbdry,\hght*\ybdry);
    
    \draw[-,black,thick](\value{co}*.5*\xbdry,\value{co}*\ybdry)--(\value{co}*.5*\xbdry+\hght*\xbdry,\value{co}*\ybdry);
    
    \draw[-,black,thick](\value{co}*.5*\xbdry,\value{co}*\ybdry)--(\value{co}*\xbdry,0);
    
    \stepcounter{co};
}
\node at (.5*\xbdry+0*\xbdry,1/3){\scalebox{1}{$\emptyset$}};
\node at (.5*\xbdry+1*\xbdry,1/3){\scalebox{.4}{\yng(2)}};
\node at (.5*\xbdry+2*\xbdry,1/3){\scalebox{.25}{\yng(4,2)}};
\node at (.5*\xbdry+3*\xbdry,1/3){\scalebox{.18}{\yng(6,4,2)}};
\node at (\xbdry+0*\xbdry,2/3){\scalebox{.4}{\yng(1)}};
\node at (\xbdry+1*\xbdry,2/3){\scalebox{.25}{\yng(3,1)}};
\node at (\xbdry+2*\xbdry,2/3){\scalebox{.18}{\yng(5,3,1)}};
\node at (\xbdry+3*\xbdry,2/3){\scalebox{.12}{\yng(7,5,3,1)}};

\node at (\xbdry+0*\xbdry,4/3){\scalebox{.4}{\yng(1,1)}};
\node at (\xbdry+1*\xbdry,4/3){\scalebox{.25}{\yng(3,1,1)}};
\node at (\xbdry+2*\xbdry,4/3){\scalebox{.18}{\yng(5,3,1,1)}};
\node at (1.5*\xbdry+0*\xbdry,5/3){\scalebox{.25}{\yng(2,1,1)}};
\node at (1.5*\xbdry+1*\xbdry,5/3){\scalebox{.18}{\yng(4,2,1,1)}};
\node at (1.5*\xbdry+2*\xbdry,5/3){\scalebox{.12}{\yng(6,4,2,1,1)}};
\node at (1.5*\xbdry+0*\xbdry,7/3){\scalebox{.25}{\yng(2,2,1,1)}};
\node at (1.5*\xbdry+1*\xbdry,7/3){\scalebox{.18}{\yng(4,2,2,1,1)}};
\node at (1.5*\xbdry+2*\xbdry,7/3){\scalebox{.12}{\yng(6,4,2,2,1,1)}};
\node at (2*\xbdry+0*\xbdry,8/3){\scalebox{.18}{\yng(3,2,2,1,1)}};
\node at (2*\xbdry+1*\xbdry,8/3){\scalebox{.12}{\yng(5,3,2,2,1,1)}};
\node at (2*\xbdry+0*\xbdry,10/3){\scalebox{.18}{\yng(3,3,2,2,1,1)}};
\node at (2.5*\xbdry+0*\xbdry,11/3){\scalebox{.12}{\yng(4,3,2,2,1,1)}};
\end{tikzpicture}}
%
%
%
%
\resizebox{.32\textwidth}{!}{
\begin{tikzpicture}[scale = 1,bull/.style={circle,draw=black,fill=black,thick,scale=.6}]
\tikzmath{\xpoint=1.732050808;\ypoint=1;\xbdry=1.154700538;\ybdry=1;\hght=8;\smallx=.8660254038;}
\clip (.8660254038,.5) circle (3.1);

\setcounter{li}{-2}
\whiledo{\value{li}<3}
{
    \setcounter{co}{-3}
    \whiledo{\value{co}<4}
    {
        \setcounter{xcurr}{2*\value{li}};
        \setcounter{ycurr}{2*\value{co}};
        \node at (\value{xcurr}*\xpoint,\value{ycurr}*\ypoint)[bull]{};
        \draw[-,black,very thin](\value{xcurr}*\xpoint+\xbdry,\value{ycurr}*\ypoint)--(\value{xcurr}*\xpoint+.5*\xbdry,\value{ycurr}*\ypoint+\ybdry)--(\value{xcurr}*\xpoint-.5*\xbdry,\value{ycurr}*\ypoint+\ybdry)--(\value{xcurr}*\xpoint-\xbdry,\value{ycurr}*\ypoint)--(\value{xcurr}*\xpoint-.5*\xbdry,\value{ycurr}*\ypoint-\ybdry)--(\value{xcurr}*\xpoint+.5*\xbdry,\value{ycurr}*\ypoint-\ybdry)--cycle;
        \stepcounter{co};
    }
    \stepcounter{li};
}

\setcounter{li}{-2};
\whiledo{\value{li}<2}
{
    \setcounter{co}{-4};
    \whiledo{\value{co}<4}
    {
        \setcounter{xcurr}{2*\value{li}+1};
        \setcounter{ycurr}{2*\value{co}+1};
        \node at (\value{xcurr}*\xpoint,\value{ycurr}*\ypoint)[bull]{};
        \draw[-,black,very thin](\value{xcurr}*\xpoint+\xbdry,\value{ycurr}*\ypoint)--(\value{xcurr}*\xpoint+.5*\xbdry,\value{ycurr}*\ypoint+\ybdry)--(\value{xcurr}*\xpoint-.5*\xbdry,\value{ycurr}*\ypoint+\ybdry)--(\value{xcurr}*\xpoint-\xbdry,\value{ycurr}*\ypoint)--(\value{xcurr}*\xpoint-.5*\xbdry,\value{ycurr}*\ypoint-\ybdry)--(\value{xcurr}*\xpoint+.5*\xbdry,\value{ycurr}*\ypoint-\ybdry)--cycle;
        \stepcounter{co};
    }
    \stepcounter{li};
}
\setcounter{co}{0}
\whiledo{\value{co}<4}
{
    \draw[-,black,thick](1.5*\value{co}*\xbdry,\value{co}*3)--(1.5*\value{co}*\xbdry+6*\xbdry,\value{co}*3+.5*\hght);
    
    \draw[-,black,thick](\value{co}*\xpoint,3*\value{co})--(\value{co}*\xpoint,\value{co});
    
    \draw[-,black,thick](\value{co}*\smallx,1.5*\value{co})--(\value{co}*\xpoint,\value{co});
    \stepcounter{co};
}
\setcounter{co}{0}
\whiledo{\value{co}<6}
{
    \draw[-,black,thick](\value{co}*\xpoint,\value{co})--(\value{co}*\xbdry+.5*\hght*\xbdry,\hght);
    \draw[-,black,thick](\value{co}*.5*\xbdry,\value{co})--(1.5*\value{co}*\xbdry,\value{co});
    
    \draw[-,black,thick](\value{co}*.5*\xbdry,\value{co})--(\value{co}*\smallx,.5*\value{co});
    \stepcounter{co};

}
\node at (.5/3*\xbdry+\smallx/3,1.5/3){\scalebox{.9}{$\emptyset$}};
\node at ({(2*\xbdry+\smallx)/3},2.5/3){\scalebox{.35}{\yng(1)}};
\node at ({(2*\xbdry+\smallx)/3},3.5/3){\scalebox{.23}{\yng(2)}};
\node at ({(\xbdry+\smallx+2*\smallx)/3},4.5/3){\scalebox{.23}{\yng(1,1)}};
\node at ({(3*\xbdry+\xbdry)/3},5/3+.1){\scalebox{.2}{\yng(2,1,1)}};
\node at ({(5*\xbdry)/3+.04},5/3+.1){\scalebox{.2}{\yng(3,1)}};
\node at ({(3.5*\xbdry+3*\smallx)/3+.1},4.5/3+.03){\scalebox{.2}{\yng(3,1,1)}};
\node at ({(3*\xbdry+\xbdry)/3+.04},7/3-.04){\scalebox{.18}{\yng(2,2,1,1)}};
\node at ({(4*\xbdry+\xbdry)/3+.02},7/3-.11){\scalebox{.16}{\yng(4,2)}};
\node at ({(3.5*\xbdry+3*\smallx)/3+.12},7.5/3-.1){\scalebox{.13}{\yng(5,3,1)}};
\node at ({(5*\xbdry+3*\smallx)/3},5.5/3-.06){\scalebox{.16}{\yng(4,2,1,1)}};
\node at ({(5*\xbdry+3*\smallx)/3-.05},6.5/3){\scalebox{.13}{\yng(5,3,1,1)}};
\node at ({(5.5*\xbdry+3*\smallx)/3-.07},7.5/3+.06){\scalebox{.11}{\yng(4,2,2,1,1)}};
\node at ({(4*\xbdry+3*\smallx)/3+.1},8.5/3-.04){\scalebox{.12}{\yng(3,2,2,1,1)}};
\node at ({(4*\xbdry+3*\smallx)/3+.1},9.5/3+.02){\scalebox{.12}{\yng(3,3,2,2,1,1)}};
\end{tikzpicture}}
\resizebox{.32\textwidth}{!}{\begin{tikzpicture}[scale = 1,bull/.style={circle,draw=black,fill=black,thick,scale=.6}]
\tikzmath{\xpoint=1;\ypoint=1;\xbdry=1;\ybdry=1;\hght=8;\smallx=1; \smally=1;}
\clip (.8660254038,.5) circle (3.1);

\setcounter{li}{-2}
\whiledo{\value{li}<3}
{
    \setcounter{co}{-3}
    \whiledo{\value{co}<4}
    {
        \setcounter{xcurr}{2*\value{li}};
        \setcounter{ycurr}{2*\value{co}};
        \node at (\value{xcurr}*\xpoint,\value{ycurr}*\ypoint)[bull]{};
        \draw[-,black,very thin](\value{xcurr}*\xpoint+\xbdry,\value{ycurr}*\ypoint+\ybdry)--(\value{xcurr}*\xpoint+\xbdry,\value{ycurr}*\ypoint-\ybdry)--(\value{xcurr}*\xpoint-\xbdry,\value{ycurr}*\ypoint-\ybdry)--(\value{xcurr}*\xpoint-\xbdry,\value{ycurr}*\ypoint+\ybdry)--cycle;
        \stepcounter{co};
    }
    \stepcounter{li};
}

\setcounter{li}{0}
\whiledo{\value{li}<3}
{
    \setcounter{co}{-1}
    \whiledo{\value{co}<\value{li}}
    {
        \setcounter{xcurr}{2*\value{li}};
        \setcounter{ycurr}{2*\value{co}+2};
        
        \draw[-,black,thick] (\value{xcurr}*\xpoint+1*\xpoint,\value{ycurr}*\ypoint+0*\xpoint) -- (\value{xcurr}*\xpoint+1*\xpoint,\value{ycurr}*\ypoint+1*\ypoint);
        \draw[-,black,thick] (\value{xcurr}*\xpoint+2*\xpoint,\value{ycurr}*\ypoint+0*\xpoint) -- (\value{xcurr}*\xpoint+2*\xpoint,\value{ycurr}*\ypoint+2*\ypoint);
        \draw[black,thick] (\value{xcurr}*\xpoint+3*\xpoint,\value{ycurr}*\ypoint+1*\xpoint) -- (\value{xcurr}*\xpoint+3*\xpoint,\value{ycurr}*\ypoint+2*\ypoint);
        
        \draw[-,black,thick] (\value{xcurr}*\xpoint+1*\xpoint,\value{ycurr}*\ypoint+1*\xpoint) -- (\value{xcurr}*\xpoint+2*\xpoint,\value{ycurr}*\ypoint+0*\ypoint);
        \draw[-,black,thick] (\value{xcurr}*\xpoint+2*\xpoint,\value{ycurr}*\ypoint+2*\xpoint) -- (\value{xcurr}*\xpoint+3*\xpoint,\value{ycurr}*\ypoint+1*\ypoint);

        \draw[-,black,thick] (\value{xcurr}*\xpoint+0*\xpoint,\value{ycurr}*\ypoint+0*\xpoint) -- (\value{xcurr}*\xpoint+2*\xpoint,\value{ycurr}*\ypoint+0*\ypoint);
        \draw[-,black,thick] (\value{xcurr}*\xpoint+1*\xpoint,\value{ycurr}*\ypoint+1*\xpoint) -- (\value{xcurr}*\xpoint+3*\xpoint,\value{ycurr}*\ypoint+1*\ypoint);
        
        \draw[-,black,thick] (\value{xcurr}*\xpoint+0*\xpoint,\value{ycurr}*\ypoint+0*\xpoint) -- (\value{xcurr}*\xpoint+2*\xpoint,\value{ycurr}*\ypoint+2*\ypoint);
        
        \stepcounter{co};
    }
    \stepcounter{li};
}

\node at (0*\xbdry+2*\smallx/3,0*\ybdry+1*\smally/3){\scalebox{.9}{$\emptyset$}};
\node at (1*\xbdry+1*\smallx/3,0*\ybdry+1*\smally/3){\scalebox{.35}{\yng(1)}};
\node at (1*\xbdry+2*\smallx/3,0*\ybdry+2.25*\smally/3){\scalebox{.3}{\yng(2,1)}};
\node at (2*\xbdry+1*\smallx/3,0*\ybdry+2*\smally/3){\scalebox{.23}{\yng(3,1,1)}};
\node at (2*\xbdry+2.25*\smallx/3,0*\ybdry+0.75*\smally/3){\scalebox{.18}{\yng(4,1,1,1)}};
\node at (1*\xbdry+2*\smallx/3,1*\ybdry+1*\smally/3){\scalebox{.3}{\yng(2,2)}};
\node at (2*\xbdry+0.75*\smallx/3,1*\ybdry+0.75*\smally/3){\scalebox{.23}{\yng(3,2,1)}};
\node at (2*\xbdry+2*\smallx/3,1*\ybdry+2.25*\smally/3){\scalebox{.18}{\yng(4,3,2,1)}};
\node at (3*\xbdry+1*\smallx/3,0*\ybdry+0.75*\smally/3){\scalebox{.15}{\yng(5,2,1,1,1)}};
\node at (3*\xbdry+2.25*\smallx/3,0*\ybdry+2.25*\smally/3){\scalebox{.13}{\yng(6,3,2,1,1,1)}};
\node at (3*\xbdry+2.25*\smallx/3,1*\ybdry+0.75*\smally/3){\scalebox{.13}{\yng(6,3,3,1,1,1)}};
\node at (3*\xbdry+1*\smallx/3,1*\ybdry+2*\smally/3){\scalebox{.18}{\yng(5,3,3,1,1)}};
\node at (2*\xbdry+2.25*\smallx/3,2*\ybdry+0.75*\smally/3){\scalebox{.18}{\yng(4,4,2,2)}};
\node at (3*\xbdry+1*\smallx/3,2*\ybdry+0.75*\smally/3){\scalebox{.15}{\yng(5,4,3,2,1)}};
\end{tikzpicture}}


\end{center}
\caption{$3$-cores in types $A_2$ and $G_2$, and self-conjugate $4$-cores in type $C_2$.}
\label{fig:a2}
\end{figure}
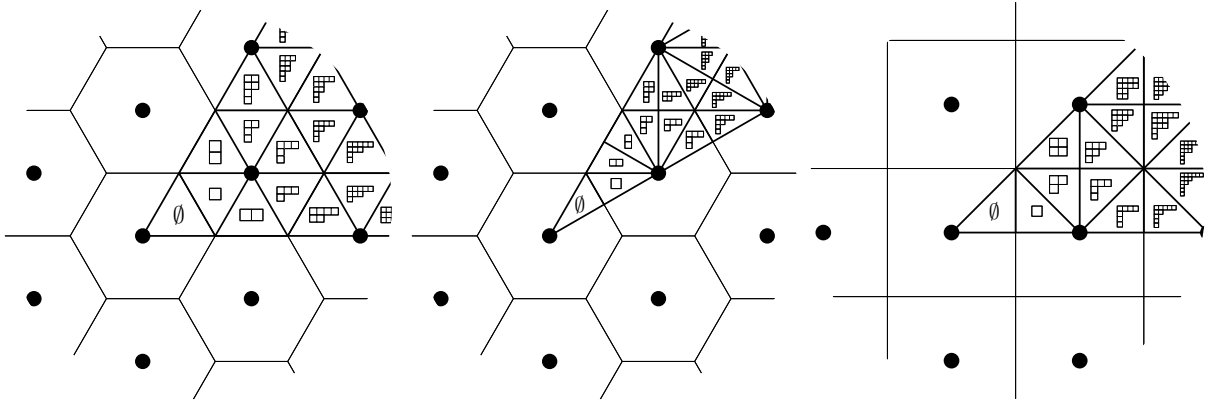

To produce similar combinatorial models for the quotients $ \widetilde{W}/W$ of other classical types ($X_n \in \{A_n,B_n,C_n,D_n\}$), we may embed $\chk{\Q}_{X_n}$ into an appropriate type $A$ coroot lattice. ~\Cref{fig:a2} illustrates these in rank 2, as well as a similar model for $X_n=G_2$.



Under the correspondence between $a$-cores and $\chk{\Q}_{a}$ of~\Cref{eq:core_to_coroot}, the set of $(a,b)$-cores turn out to be exactly those coroot points that sit inside of a certain affine transformation of the fundamental alcove $\Sommers(b)$, which includes a $b$-fold dilation, called the $b$-Sommers region (see \Cref{def:sommers}).  The natural generalization of $\core(a,b)$ to any affine Weyl group is therefore the intersection of the coroot lattice $\chk{\Q}_{X_n}$ with its $b$-Sommers region, so that $\core(a,b) = \core(A_{a-1},b)$. In other words,
\begin{equation}
\label{eq:Wcores}
    \core(X_n,b) := \chk{\Q}_{X_n} \cap \Sommers_{X_n}(b).
\end{equation}

\subsection{Previous Work}

Under the bijections of~\Cref{eq:core_to_coroot}, we noticed in~\cite{thiel2017strange} that the number of boxes in $\lambda$ could be computed from the coroot $q_\lambda$ as described above, or the inversion set of $\widetilde{w}_\lambda^{-1}$, where $\mathsf{inv}(\widetilde{w})=\widetilde{\Phi}^+ \cap \widetilde{w}(-\widetilde{\Phi}^+)$. More precisely:

\begin{proposition}[{\cite[Proposition 6.4 \& Corollary 6.7]{thiel2017strange}}] Let $\lambda$ be an $a$-core and $\chk{\rho}$ be the sum of the fundamental coweights in type $A_{a-1}$. Then
$$\size(\lambda) = \sum\limits_{\alpha+k\delta\in\mathsf{inv}(\waf^{-1}_\lambda)}k = \left\langle \frac{a}{2}q_{\lambda} - \chk{\rho},q_\lambda \right\rangle.$$
\end{proposition}
It was natural to consider the corresponding statistic in \emph{any} affine Weyl group $\widetilde{W}(X_n)$ acting on $V$, restricting to a certain finite set of coroots $\core(X_n,b)$ (defined below in~\Cref{eq:Wcores}, in analogy with simultaneous $(a,b)$-cores).  The latter two authors showed that for simply-laced Weyl groups, the result mirrored~\Cref{thm:arm}.

\begin{theorem}[{\cite[Theorem 1.10]{thiel2017strange}}]
\label{thm:us}
Let $X_n$ be a simply-laced Cartan type with Coxeter number $h$, and let $b$ be coprime to $h$.  Then
\begin{equation*} \Expt{q \in \core(X_n,b)}{\size(q)} = \frac{n(b-1)(h+b+1)}{24}. \end{equation*}
\end{theorem}

When applied to $X_n=A_{a-1}$ (so that $n=a-1$ and $h=a$), this result gives a proof of the left equality of \Cref{thm:arm} for the expected size of simultaneous $(a,b)$-cores.  But since self-conjugate cores are a combinatorial model for coroots in the \emph{non-simply-laced} type $C_n$, we were unable to similarly specialize \Cref{thm:us} to conclude the right equality of \Cref{thm:arm} for the expected size of a self-conjugate simultaneous core.

\subsection{Improved Size Statistic}
In this paper, we describe a modification of the $\size$ statistic to incorporate the lengths of the roots. This appears advantageous over the original statistic of \cite{thiel2017strange}; we are able to apply the Ehrhart-theoretic techniques of P.~Johnson outside of simply-laced type. Normalize root systems so that the highest root has length $2$, and write $r$ for the ratio of the length of a long to a short root. For $\waf \in \widetilde{W}/W$, define 
\begin{equation}
\label{eqn:size_def}
    \chk{\size}(\waf) := \left(\sum_{\substack{\alpha+k\delta \in \mathsf{inv}(\waf^{-1})\\\alpha \text{ long}}}k\right)+r\left(\sum_{\substack{\alpha+k\delta \in \mathsf{inv}(\waf^{-1})\\\alpha \text{ short}}}k\right) 
\end{equation}

This recovers the original statistic $\size$ in simply-laced type, but disagrees in non-simply-laced type when $r>1$.  A similar statistic was independently considered in~\cite{chapelier2022atomic}.

Using a bijection analogous to those of \Cref{eq:core_to_coroot}, we interpret $\chk{\size}$ as statistics on the combinatorial models of~\Cref{sec:models}. For instance, we shows that $\chk{\size}$ in type $C_n$ corresponds to the number of boxes in the corresponding self-conjugate $2n$-core (see~\Cref{fig:a2}).

Following the same strategy as in~\cite{thiel2017strange}, we find an affine Weyl group element that maps $\Sommers(b)$ to a $b$-fold dilation of the fundamental alcove (correctly modifying the $\chk{\size}$ statistic), and then apply Ehrhart theory to compute the expected value of $\chk{\size}$ on $\core(X_n,b)$.

\begin{restatable}{theorem}{mainthm}
\label{thm:main_thm}
For $X_n$ an irreducible rank $n$ Cartan type with root system $\Phi$,
\begin{align*}
\Expt{q \in \core(X_n,b)}{\chk{\size}(q)}=\frac{r \chk{g}}{h}\frac{n(b-1)(h+b+1)}{24},
\end{align*}
where $h$ is the Coxeter number of $X_n$, $\chk{g}$ is the dual Coxeter number for $\chk{\Phi}$, and $r$ is the ratio of the length of a long root to the length of a short root in $\Phi$.
\end{restatable}

The extra factor of $\frac{r\chk{g}}{h}$ is \emph{invisible} in the simply-laced case, where $\chk{\Phi}=\Phi$, $\chk{g}=h$, and $r=1$.  As an immediate application of~\Cref{thm:main_thm}, we conclude both equalities in~\Cref{thm:arm} by specializing to these types.    Interestingly, although the expected number of boxes in a simultaneous core and in a self-conjugate simultaneous core happen to be the same, the formulas have quite different interpretations: the factor of $a-1$ corresponds to the dimension $n$ for ordinary simultaneous cores, but to $\chk{g}$ in the self-conjugate case.   

We prove \Cref{thm:main_thm} for non-simply-laced types in \Cref{sec:expectation}, after some setup, including a careful definition of $\Sommers(b)$. Along the way, we briefly generalize other results from~\cite{thiel2017strange}, including in particular~\Cref{thm:max_size} concerning the maximum size.

\section{Background}

We give a brief account of most of the notation used in the remainder of the paper for objects associated to affine root systems. For definitions and greater detail, we refer the reader to standard references (e.g. \cite{humphreys1992reflection}) or to the previous paper of the second and third author \cite[Section 2]{thiel2017strange}.

\subsection{Root Systems}
\label{sec:root_systems}

Let $V$ be a Euclidean space of dimension $n$, and $\Phi\subseteq V$ be an irreducible crystallographic root system in $V$ of type $X_n$.  We often suppress the $X_n$ notation when there is only one root system under consideration.  Denote a system of simple roots by $\simp=\{\alpha_1,\dots, \alpha_n\}$, and the corresponding positive roots by $\Phi^+$.

For any $\alpha\in\Phi$, we may write $\alpha$ in the basis of simple roots as $\alpha=\sum_{i=1}^n a_i\alpha_i$, where the coefficients $a_i$ are either all nonnegative or all nonpositive. The \defn{height} of $\alpha$ is the sum of the coefficients:
$\h(\alpha):=\sum_{i=1}^n a_i$.
Notice that $\h(\alpha)>0$ if and only if $\alpha\in\Phi^+$ and $\h(\alpha)=1$ if and only if $\alpha\in\simp$.
There is a unique root $\amax$ of maximal height called the \defn{highest root} of $\Phi$, and we denote its coefficients by $c_i$, that is, $\amax =\sum_{i=1}^n c_i\alpha_i \in\Phi$. In addition, the \defn{Coxeter number} of $\Phi$ is $h:=1+\h(\amax)=1+\sum_{i=1}^n c_i$.

For a root $\alpha\in\Phi$, define its \defn{coroot} as $\check{\alpha}:=\frac{2\alpha}{\|\alpha\|^2}$.
Define the \defn{dual root system} of $\Phi$ as $\chk{\Phi}:=\left\{\chk{\alpha}:\alpha\in\Phi\right\}$. It is itself an irreducible crystallographic root system, and hence also has a highest root $\widetilde{\gamma}$; note that although $\widetilde{\gamma}$ is by definition the coroot of some $\alpha\in\Phi$, this $\alpha$ is typically not the highest root $\amax$. Writing $\widetilde{\gamma}=\sum_{i=1}^n d_i\chk{\alpha_i}$ as a sum of the simple coroots in $\chk{\Phi}$, then we define the \defn{dual Coxeter number} $\chk{g}:=1+\sum_{i=1}^n d_i$.

Define the \defn{coroot lattice} $\chk{\Q}$ of $\Phi$ as the lattice in $V$ generated by $\chk{\Phi}$. Finally, let $(\cw_1,\cw_2,\ldots,\cw_n)$ be the basis that is dual to the basis $(\alpha_1,\alpha_2,\ldots,\alpha_n)$ of $V$ consisting of the simple roots, so that $\langle \cw_i,\alpha_j\rangle=\delta_{i,j}$.
Then $\cw_1,\cw_2,\ldots,\cw_n$ are the \defn{fundamental coweights}. They are a basis of the \defn{coweight lattice} 
\[\chk{\Lambda}:=\{x\in V:\langle x,\alpha\rangle\in\Z\text{ for all }\alpha\in\Phi\}\]
of $\Phi$, which contains $\chk{\Q}$ as a sublattice. The sum of these basis elements, $\chk{\rho}=\sum_{i=1}^n\cw_i$, will be of particular importance. For notational convenience, we define $\chk{\omega}_0:=0$.

\begin{convention}\label{con:normalize}
We normalize the inner product $\langle\cdot,\cdot\rangle$ on $V$ so that $
\langle\amax,\amax\rangle=2$ and call $\alpha \in \Phi$ a \defn{long root} if $\langle \alpha,\alpha \rangle=2$.
\end{convention}

In particular, all long roots are their own coroots.  A \defn{short root} is a root with $\langle \alpha,\alpha \rangle<2$.  If the system has short roots $\alpha$, then $\chk{\alpha}=r\alpha$ for an integer $r\in \{2,3\}$ independent of $\alpha$. If the system does not have short roots, it is called \defn{simply-laced}.  Note that $\chk{\Phi}$ is itself a root system, but not subject to \Cref{con:normalize}.

\subsection{Affine Weyl Groups and Affine Root Systems}
\label{sec:affine}

The \defn{Weyl group} $W$ associated to a root system $\Phi$ is the subgroup of $\mathrm{GL}(V)$ generated by the \defn{simple reflections} \[s_i=s_{\alpha_i} : x\mapsto x - 2\frac{\langle \alpha_i,x\rangle}{\langle \alpha_i,\alpha_i\rangle}\alpha_i\] for $\alpha_i \in \Delta$. The corresponding \defn{affine Weyl group} $\wa$ is the subgroup of distance-preserving transformations on $V$ generated by the simple reflections $\{s_\alpha\}_{\alpha \in \Delta}$ together with the additional \defn{affine simple reflection} \[s_0: x\mapsto x - (\langle \widetilde{\alpha},x\rangle-1)\widetilde{\alpha}.\]

One readily checks that the affine Weyl group $\wa$ acts on both $\chk{\Q}$ and $\chk{\Lambda}$. For any $y\in V$, there is an associated translation $t_y:x\mapsto x+y$. If we identify $\chk{\Q}$ with the corresponding group of translations acting on $V$, then $\wa$ may be written as the semidirect product $\wa=W\ltimes\chk{\Q}$.  For $\waf\in\wa$ we will use the notation $\waf=w\cdot t_q$ to denote this semidirect product decomposition.  This decomposition gives a bijection $W\backslash \wa\to\chk{\Q}$ given by $\waf\mapsto q$, but we will make frequent use of the following more interesting bijection:

\begin{theorem}
\label{thm:coset_to_coroot}
The map $\wa = W \ltimes \chk{\Q} \to \chk{\Q}$ defined by $\waf \mapsto \waf(0)$ descends to a $\wa$-equivariant bijection on the cosets $\wa/W$.
\end{theorem}
\begin{proof}
Evidently the first map is $\wa$-equivariant, and because $g\in W$ implies $g(0)=0$, we have that $\widetilde{w}$ and $\widetilde{w}g$ have the same image. Hence we have a well-defined equivariant map on cosets, and evidently $q\mapsto t_qW$ is its inverse, as desired.
\end{proof}

The \defn{affine root system} is defined by $\widetilde{\Phi} = \Phi\times\mathbb{Z}$, and---writing $\delta$ for a formal variable to keep track of the coefficient of $\mathbb{Z}$---we use the notation $\alpha+k\delta$ for a typical element of $\widetilde{\Phi}$.  The root system $\Phi$ embeds in $\widetilde{\Phi}$ by writing $\alpha_i$ as $\alpha_i+0\cdot \delta$, and we  define $\alpha_0:=-\amax+\delta$.  The affine Weyl group $\wa$ acts on $\widetilde{\Phi}$ by
\[ \waf \cdot (\alpha + k \delta) := w(\alpha) + (k-\langle \alpha,q\rangle)\delta, \]
where $\waf = w \cdot t_q$.

\begin{definition}
Given a reduced word $\mathsf{\widetilde{w}}=s_{i_1} s_{i_2} \cdots s_{i_\ell}$ for  $\waf \in \wa$, we define its \defn{inversion sequence} 
\[\inv(\mathsf{\widetilde{w}})=\beta_1+k_1\delta,\beta_2+k_2\delta,\ldots,\beta_\ell+k_\ell\delta,\] 
where $\beta_j+k_j\delta$ are the affine roots $(s_{i_1}\cdots s_{i_{j-1}})(\alpha_{i_j})$.
\end{definition}

There may be many reduced expressions---and hence many inversion sequences---for a given $\waf\in\wa$, but these differ only by a reordering: they record the affine hyperplanes that separate $w(\mathcal{A})$ from the fundamental alcove $\mathcal{A}$. 

\section{Core Partitions and the Type A Coroot Lattice}

As discussed in the introduction, there is a close relation between the coroot lattice for type $A_n$ and certain kinds of partitions. Much of the work in this section is well-known~\cite{armstrong2014results,fishel2010bijection,thiel2017strange}, with the exception of the $\chk{size}_i$-refinement of~\Cref{prop:ip_content}.  We also refer the reader to the recent preprint~\cite{chapelier2022atomic} and to our previous FPSAC abstract on this work~\cite{stucky85strange}.

\subsection{Coroots and Cores}
\label{sec:abaci}
In type $A_{a-1}$, one choice of simple roots is $\alpha_i:=e_{i+1}-e_{i}$ for each $1\leq i < a$. Then the highest root is $\widetilde{\alpha}=e_a-e_1$, and the coroot lattice\footnote{For safety---even though roots and coroots can be identified in type $A$---we already throw in the distinguishing check.} is
$\chk{\Q}_{a} = \chk{\Q}_{A_{a-1}} := \left\{q=(q_1,q_2\ldots,q_a) \in \mathbb{Z}^a : \sum_{i=1}^a q_i = 0 \right\}.$

An integer partition $\lambda$ can be characterized by its \defn{boundary word}---a bi-infinite sequence of \defn{beads}, which are either $\bullet$s or $\circ$s, that begins with an infinite sequence of only $\bullet$s and ends with an infinite sequence of only $\circ$s. This word encodes the boundary of $\lambda$ (in English notation) by detailing the steps taken when traversing from bottom left to top right: $\bullet$s representing steps up and $\circ$s representing steps right.  For example, the boundary word for the partition on the left of~\Cref{fig:abacus} is read from south-west to north-east as $\cdots\bullet\bullet\circ\bullet\bullet\circ\circ\bullet\circ\circ\bullet\circ\circ\cdots$.

Partitioning the boundary word into consecutive subsequences of length $a$ and stacking them vertically gives the \defn{$a$-abacus} representation of $\lambda$.  This is illustrated in the middle of~\Cref{fig:abacus}.  Finally, an $a$-abacus is called \defn{balanced} if we can draw a horizontal line between two rows with as many $\circ$s above the line as $\bullet$s below; every partition has a unique representation as a balanced $a$-abacus. 

An integer partition $\lambda$ is an $a$-core if and only if its $a$-abacus representation is \defn{flush}---that is, if each of the vertical ``runners'' of the abacus consists of an infinite sequence of only $\bullet$s followed by an infinite sequence of only $\circ$s. A flush, balanced $a$-abacus encodes a \defn{coroot} as the $a$-tuple of signed distances from beneath the lowest $\bullet$ in each runner to the line witnessing the balanced condition---the balanced condition ensures that these distances sum to zero.  We will say that a bead is at \defn{level} $\ell$ if the distance from beneath the bead to the line witnessing the balanced condition is $\ell$; note that this means that levels increase when reading down the abacus. This is illustrated on the right of~\Cref{fig:abacus}.

\begin{figure}[htbp]
\[
\text{$3$-core } \lambda=\raisebox{-.5\height}{\begin{tikzpicture}[scale=1,bull/.style={circle,draw=black,fill=black,thick,scale=.5},circ/.style={circle,draw=black,fill=white,thin,scale=.5}] \node (y) at (0,0) {\young(01201,201,1,0)};
  \node at (-1.15,-1.6)[bull]{};
  \node at (-1.15,-1.15)[bull]{};
  \node at (-.9,-.95)[circ]{};
  \node at (-.45,0)[circ]{};
  \node at (0,0)[circ]{};
  \node at (.45,.45)[circ]{};
  \node at (.9,.45)[circ]{};
  \node at (1.35,.9)[circ]{};
  \node at (1.8,.9)[circ]{};
  \node at (-.70,-.70)[bull]{};
  \node at (-.70,-.25)[bull]{};
  \node at (.25,.25)[bull]{};
  \node at (1.15,.7)[bull]{};
  \end{tikzpicture}} 
\text{ abacus } ~~
\raisebox{-.5\height}{\begin{tikzpicture}[scale=1,bull/.style={circle,draw=black,fill=black,thick,scale=.6},circ/.style={circle,draw=black,fill=white,thin,scale=.6}]
\draw[->,black,thin](-.25,-2.25)--(-.25,-5.5);
\draw[->,black,thin](.25,-2.25)--(.25,-5.5);
\draw[->,black,thin](.75,-2.25)--(.75,-5.5);
\node at (-.25,-2.5)[bull]{};
\node at (.25,-2.5)[bull]{};  
\node at (.75,-2.5)[bull]{};

\node at (-.25,-3)[bull]{};
\node at (.25,-3)[bull]{};  
\node at (.75,-3)[circ]{};

\node at (-.25,-3.5)[bull]{};
\node at (.25,-3.5)[bull]{};  
\node at (.75,-3.5)[circ]{};

\node at (-.25,-4)[circ]{};
\node at (.25,-4)[bull]{};  
\node at (.75,-4)[circ]{};

\node at (-.25,-4.5)[circ]{};
\node at (.25,-4.5)[bull]{};  
\node at (.75,-4.5)[circ]{};

\node at (-.25,-5)[circ]{};
\node at (.25,-5)[circ]{};  
\node at (.75,-5)[circ]{};
\draw[-,blue,thick](-.5,-3.75)--(1,-3.75);
\end{tikzpicture}}
\quad \text{ and coroot }
q=(0,2,-2)=-\alpha_2.\]
\caption{An example of the bijection between $a$-cores, abaci, and ${\protect \chk{\Q}_{A_{a-1}}}$ (for $a=3$).}
\label{fig:abacus}
\end{figure}

By the discussion above, $\chk{\Q}_{a}$ is in bijection with the set of $a$-cores $\core(a)$.

\begin{definition}\label{def:bij}
For $q \in \chk{\Q}_a$, we write $\lambda_q$ for the $a$-core obtained by building the flush, balanced $a$-abacus with levels of the lowest $\bullet$ in each runner given by the coordinates of $q$, and then reading this as the boundary word of a partition; for $\lambda \in \core(a)$, we write $q_\lambda$ for the corresponding coroot in $\chk{\Q}_a$ obtained by reading the boundary word of $\lambda$, producing the corresponding $a$-core, and then reading off the levels of the lowest $\bullet$ in each runner.
\end{definition}

The action of the affine symmetric group $\widetilde{\mathfrak{S}}_a=\wa(A_{a-1})$ on $\chk{\Q}_a$ is generated by the usual simple reflections $s_i$ interchanging the $i^\text{th}$ and $(i+1)^\text{st}$ positions, along with the additional affine simple reflection $s_0$: 
\begin{align*}
    s_i(q_1,\ldots,q_i,q_{i+1},\ldots,q_a) &= (q_1,\ldots,q_{i+1},q_{i},\ldots,q_a), \text{ and}\\
    s_0(q_1,\ldots,q_a) &= (q_a+1,\ldots,q_1-1).
\end{align*}  

We can translate this action of $\widetilde{\mathfrak{S}}_a$ to the set of $a$-cores~\cite[Section 2.7]{james1981representation}~\cite{lascoux2001ordering}. We think of a partition as an order ideal in $\mathbb{N} \times \mathbb{N}$ (top-left justified), where each $(i,j) \in \mathbb{N} \times \mathbb{N}$ is indexed by its \defn{content} $(i-j) \mod a$.   For $0 \leq i < a$, let the simple reflection $s_i$ act on a partition by toggling all possible boxes with content $i$ mod $n$---that is, adding all possible missing boxes with content $i$ which produce a valid Young diagram, or removing all possible present boxes with content $i$ which produce a valid Young diagram.  This extends to an action of the full affine symmetric group $\widetilde{\mathfrak{S}}_a$ on $a$-cores.

\begin{theorem}
\label{thm:sa_eq}
The action of the affine symmetric group $\widetilde{\mathfrak{S}}_a$ is preserved under the bijection between $\chk{\Q}_a$ and $\core(a)$  of~\Cref{def:bij}.  That is, for $0 \leq i <a$, $q \in \chk{\Q}_a$, and $\lambda \in \core(a)$, we have \[s_i(q) = s_i(\lambda_q) \text{ and } s_i(\lambda) = s_i(q_\lambda).\]
\end{theorem}

\subsection{Two Size Statistics}
\label{sec:size}

For $\lambda$ a partition, write: $\lambda^\intercal$ for its conjugate; $\size_i(\lambda)$ for the the number of boxes in $\lambda$ with content $i \mod a$; and $\size(\lambda)$ for the total number of its boxes.  Under the bijection between coroots and $a$-cores, we can interpret these definitions in the language of the coroot lattice. For $q=(q_1,\ldots,q_a) \in \chk{\Q}_a$, write $q^\intercal:=(-q_a,\ldots,-q_1)$ and define \begin{align*}
    \chk{\size}_i(q):=\left\langle \frac{1}{2} q-\chk{\omega}_i,q\right\rangle \text{ and }    \chk{\size}(q):=\sum_{i=1}^{a-1} \chk{\size}_i(q)=\left\langle \frac{a}{2} q-\chk{\rho},q\right\rangle.
\end{align*}

Recall that in type $A_{a-1}$ (up to the usual normalization that the sum of the entries ought to be zero), we have for $0 \leq i < a$:
\begin{align*}
    \chk{\omega}_i &= \sum_{j=1}^i e_i = (\underbrace{1,1,\ldots,1}_{i \text{ ones}},0,0,\ldots,0) \\
    \chk{\rho} &= \sum_{i=1}^{a-1} \chk{\omega}_i = (a-1,a-2,\ldots,1,0).
\end{align*}
Recall that we previously defined $\omega_0=0$, which agrees with this description.  Note that we are able to safely ignore the normalization on $\chk{\omega}_i$ and $\chk{\rho}$ because it is still enforced on the coroot $q$ when computing $\chk{\size}_i$. For instance,
\[\left\langle \frac{1}{2} q-(\chk{\omega}_i+t(1,\dots,1)),q\right\rangle = \chk{\size}_i(q) - t\langle(1,\dots, 1), q\rangle = \chk{\size}_i(q).\]

\def\numz{{\textcolor{blue}{0}}}
\def\numo{{\textcolor{green}{1}}}
\def\numt{{\textcolor{red}{2}}}
\begin{example}
\label{ex:a2size}
Continuing the example from \Cref{fig:abacus}, the coroot $q=(0,2,-2)\in \chk{\Q}_3$ corresponds to the $3$-core $\lambda=\raisebox{\height}{\young(\numz\numo\numt\numz\numo,\numt\numz\numo,\numo,\numz)}$.  This $\lambda$ has four boxes with content $\numz \mod 3$, four with content $\numo \mod 3$, and two with content $\numt \mod 3$.  We compute

\begin{align*}
    \chk{\size}_\numz(q) &= \left\langle \frac{1}{2}(0,2,-2),(0,2,-2)\right\rangle = 4 = \size_\numz(\lambda),\\
    \chk{\size}_\numo(q) &= \left\langle \frac{1}{2}(0,2,-2)-\textcolor{green}{(1,0,0)},(0,2,-2)\right\rangle = 4=\size_\numo(\lambda),\\
     \chk{\size}_\numt(q) &= \left\langle \frac{1}{2}(0,2,-2)-\textcolor{red}{(1,1,0)},(0,2,-2)\right\rangle = 2=\size_\numt(\lambda).
\end{align*}
\end{example}

This correspondence holds in general, as follows.

\begin{proposition}
\label{prop:ip_content}
For any $a$-core $\lambda$ and $q=q_\lambda$, then $\lambda^\intercal = \lambda_{q^\intercal}$ and $\chk{\size}(\lambda)=\size(q)$.  In fact, for any $0\leq i< a$, $\size_i(\lambda)=\chk{\size}_i(q)$.
\end{proposition}
\begin{proof}
The statement about conjugation follows by observing that the boundary word of a partition and its conjugate are related by reversing and interchanging $\bullet \leftrightarrow \circ$.

Write $q=(q_1,\dots, q_a)$ and $\lambda=\lambda_q$. We compute directly that 
\begin{align*}
\left\langle \frac{a}{2} q-\chk{\omega}_i, q \right\rangle
&= \sum_{j=1}^i \left(\frac{1}{2}q_j-1\right)q_j + \sum_{j=i+1}^a \frac{1}{2}q_j^2 \\
&= \sum_{j=1}^i \left(\frac{(q_j-1)q_j}{2} - \frac{q_j}{2}\right) + \sum_{j=i+1}^a \left(\frac{(q_j+1)q_j}{2} - \frac{q_j}{2}\right) \\
&= \sum_{j=1}^i \frac{q_j(q_j-1)}{2} + \sum_{j=i+1}^a \frac{(q_j+1)q_j}{2}.
\end{align*}
Observing that this is a sum over runners of certain triangular numbers, it would suffice to show the boxes of content $i$ in $\lambda$ may be partitioned so each bead on runner $j$ at level $\ell$ corresponds to:
$$
\begin{cases}
    \text{$\ell$ boxes} & \text{$\ell>0$, black bead, $j\leq i$}, \\
    \text{$-\ell$ boxes} & \text{$\ell\leq 0$, white bead, $j\leq i$}, \\
    \text{$\ell-1$ boxes} & \text{$\ell>0$, black bead, $j> i$}, \\
    \text{$-\ell+1$ boxes} & \text{$\ell\leq 0$, white bead, $j> i$}.
\end{cases}
$$

To do this, begin by partitioning the diagram beneath the main diagonal; that is, into the boxes with negative content and non-negative content.   A black bead on runner $j$ at level $\ell>0$ represents a vertical edge on the boundary of $\lambda$, and every box not beneath the main diagonal is in the same row as some such edge. In that row, on or above the main diagonal, there is one box each of content $0,1,2,\dots,(\ell-1)a+j-1$. Counting the number of such boxes with content $i \bmod a$, we see there are $\ell-1$ of them if $j\leq i$, or $\ell$ if $j> i$.


Similarly, a white bead on runner $j$ at level $\ell\leq 0$ represents a horizontal edge on the boundary of $\lambda$, and every box beneath the main diagonal is in the same row as some such edge. In that row, beneath the main diagonal, there is one box each of content $-1,-2,\dots,(\ell-1)a+j$. Counting the number of such boxes with content $i \bmod a$, we see there are $-\ell$ of them if $j\leq i$, or $-\ell+1$ if $j>i$.
\end{proof}

\section{Size Statistics in General Type}
\label{sec:size_gen}

We now turn to the general definition of the size statistic for affine Weyl groups---note that when we leave type $A$, we do not have a uniform combinatorial interpretation (although see the next \Cref{sec:models} for interpretations in the other classical types $B_n,C_n,D_n$ and in type $G_2$).

\begin{definition}
\label{def:sizewords}
Fix $\waf \in \wa$ and a reduced word $\mathsf{\widetilde{w}}^{-1}=a_1 \cdots a_{\ell}$ for $\waf^{-1}$, with inversion sequence $\inv(\mathsf{\widetilde{w}}^{-1})=\beta_1+k_1\delta,\beta_2+k_2\delta,\ldots,\beta_\ell+k_\ell\delta.$   For any $i \in \{0,1,\ldots,n\}$ with corresponding simple reflection $s_i$ and simple root $\alpha_i$, define \begin{align*}  \chk{\size}_i(\mathsf{\widetilde{w}})&=\frac{2}{\left\langle \alpha_i, \alpha_i\right\rangle}  \sum_{\substack{1\leq j \leq \ell \\ a_j=s_i}} k_j .\end{align*}
\end{definition}

\begin{example}
\label{ex:a2}
Continuing~\Cref{ex:a2size}, the coroot $q=(0,2,-2)=-2\chk{\alpha_2}$ corresponds to the the coset containing $\waf = s_1s_0s_1s_2s_1s_0 \in \widetilde{A}_2$. We compute the inversion sequence for the reduced word $\mathsf{\widetilde{w}^{-1}} = s_{\textcolor{blue}{0}}s_{\textcolor{green}{1}}s_{\textcolor{red}{2}}s_{\textcolor{green}{1}}s_{\textcolor{blue}{0}}s_{\textcolor{green}{1}}$ representing $\waf^{-1}$:
\[
-\widetilde{\alpha}+\textcolor{blue}{1}\cdot\delta,~~
-\alpha_2+\textcolor{green}{1}\cdot\delta,~~
-\widetilde{\alpha}+\textcolor{red}{2}\cdot\delta,~~
-\alpha_1+\textcolor{green}{1}\cdot\delta,~~
-\widetilde{\alpha}+\textcolor{blue}{3}\cdot\delta,~~
-\alpha_2+\textcolor{green}{2}\cdot\delta.\]
We observe that $\size_{\textcolor{blue}{0}}(\mathsf{\widetilde{w}})=\textcolor{blue}{1}+\textcolor{blue}{3}=4$, 
    $\size_{\textcolor{green}{1}}(\mathsf{\widetilde{w}})=\textcolor{green}{1}+\textcolor{green}{1}+\textcolor{green}{2}=4$, and 
    $\size_{\textcolor{red}{2}}(\mathsf{\widetilde{w}})=\textcolor{red}{2},$
agreeing with the previously-computed $\size_i(q)$ and $\size_i(\lambda_q)$.
\end{example}

\Cref{def:sizewords} turns out to not depend on the reduced word chosen for $\waf^{-1}$, or on the choice of coset representative.

\begin{proposition}
\label{prop:size_welldef}
Let $\widetilde{w}$, $\widetilde{w}'$ represent the same coset of $\widetilde{W}$, and let $\mathsf{\widetilde{w}}$ and $\mathsf{\widetilde{w}}'$ be any two reduced words for those elements.  Then $\chk{\size}_i(\mathsf{\widetilde{w}})=\chk{\size}_i(\mathsf{\widetilde{w}}')$.
\end{proposition}

\begin{proof}
We first show that $\chk{\size}_i$ is constant when $\widetilde{w}=\widetilde{w}'$.  Since the set of reduced words of $\waf$ are connected under braid moves, it suffices to show this when $\mathsf{\widetilde{w}}$ and $\mathsf{\widetilde{w}}'$ differ by a single braid move.  In that case, we have one of
\begin{enumerate}[(i)]
\item $\mathsf{\widetilde{w}} = \cdots (s_i s_j s_i) \cdots$ and $\mathsf{\widetilde{w}}' = \cdots (s_j s_i s_j) \cdots$, 
\item $\mathsf{\widetilde{w}} = \cdots (s_i s_j s_i s_j) \cdots$ and $\mathsf{\widetilde{w}}' = \cdots (s_j s_i s_j s_i) \cdots$, or
\item $\mathsf{\widetilde{w}} = \cdots (s_i s_j s_i s_j s_i s_j) \cdots$ and $\mathsf{\widetilde{w}}' = \cdots (s_j s_i s_j s_i s_j s_i) \cdots$,
\end{enumerate}
corresponding to a braid move of type $A_2, B_2,$ or $G_2$.  In each case, the order of the corresponding roots in the rank two parabolic subgroup is reversed; since the positions of the $s_i$ are also reversed in cases (ii) and (iii), these are immediate.  And in case (i), the statement follows because the roots in the rank two parabolic (of type $A_2$) are of the form $\alpha,\alpha+\beta,\beta$.

Finally, if $\waf = \widetilde{v}s_i$ for $i\neq 0$, then the affine roots $\inv(\waf^{-1})$ are simply $\alpha_i$ and $s_i(\beta_i)+k_i\delta$ for $\beta_i+k_i\delta\in\inv(\widetilde{v}^{-1})$. By induction, the size of $\waf$ is invariant under right-multiplication by $W$-elements, as needed.
\end{proof}

For any coset representative $\waf$, we may therefore define the statistic $\chk{\size}$ on $\wa/W$ as
\[\chk{\size}(\waf) := \sum_{i=0}^n \chk{\size}_i(\waf) =  \left(\sum_{\substack{\alpha+k\delta \in \mathsf{inv}(\waf^{-1})\\\alpha \text{ long}}}k\right)+r\left(\sum_{\substack{\alpha+k\delta \in \mathsf{inv}(\waf^{-1})\\\alpha \text{ short}}}k\right).\]

\begin{definition}
\label{def:sizelattice}
Recall that we expand the highest root as a sum of simples as $\widetilde{\alpha}=\sum_{i=1}^n c_i \alpha_i$, and set $c_0:=1$.  For $q \in \chk{\Q}$, define
\[\chk{\size_i}(q)=\Big\langle \frac{c_i}{2}q-\chk{\omega_i},q\Big\rangle.\]
\end{definition}

Because $\frac{2}{\left\langle \alpha_i, \alpha_i\right\rangle} = 1$ if $\alpha_i$ is long, and is $r$ if $\alpha_i$ is short, and because $\sum_{i=0}^n c_i = h$ and $\sum_{i=0}^n \chk{\omega}_i = \chk{\rho}$, we obtain 
\[\chk{\size}(q) = \sum_{i=0}^n \chk{\size}_i(q) = \left\langle \frac{h}{2}q - \chk{\rho},q \right\rangle.\]

\begin{theorem}
\label{thm:sizer}
For $\waf = w t_q \in \wa= W \ltimes \chk{\Q}$, we have $\chk{\size}_i(\waf) =\chk{\size_i}(w(q))$ and $\chk{\size}(\waf)=\chk{\size}(w(q)).$
\end{theorem}

Note that $w(q) = wt_q(0)$, and so this theorem states that $\chk{\size}_i$ and $\chk{\size}$ are preserved under the equivariant bijection defined in \Cref{thm:coset_to_coroot}.

\begin{proof}
It suffices to prove the statement for $\chk{\size}_i$.  Let $j \neq 0$ and let $i \in \{0,1,\ldots,n\}$.  We compute $\chk{\size}_i(s_jw(q))$:
\begin{align*}
    \chk{\size}_i(s_jw(q))
    &= \left\langle \frac{c_i}{2} s_jw(q)-\chk{\omega}_i,s_jw(q)\right\rangle \\
    &= \left\langle \frac{c_i}{2} \left[w(q)-\langle \chk{\alpha}_j,w(q)\rangle\alpha_j\right]-\chk{\omega}_i,\left[w(q)-\langle \chk{\alpha}_j,w(q)\rangle\alpha_j\right]\right\rangle \\
    &=\left\langle \frac{c_i}{2} q-\chk{\omega}_i,w(q)\right\rangle+\left\langle \chk{\alpha}_j, w(q)\right\rangle \cdot \left\langle \chk{\omega}_i,\alpha_j\right\rangle\\
    &=\chk{\size}_i(w(q))+\begin{cases}\big\langle \chk{\alpha}_j, w(q)\big\rangle & \text{if } i=j \\ 0 & \text{if } i\neq j\end{cases}.
\end{align*}
    Similarly, we compute $\chk{\size}_i(s_0w(q)) = \chk{\size}_i(w(q))+\begin{cases}1-\left\langle \widetilde{\alpha}, w(q)\right\rangle & \text{if } i=0 \\ 0 & \text{if } i\neq 0\end{cases}$.

We now argue by induction on the length of an affine element $\waf$, with base case coming from the identity $e \leftrightarrow 0$ giving $\chk{\size}_i(e)=\chk{\size}_i(0)=0$.  Consider now an affine element $\waf=w\cdot t_{-q}$ of length $\ell-1$, a reduced expression $\waf=a_1a_2\cdots a_{\ell-1}=w \cdot t_q$, and simple transposition $s_j$. 

The result follows by comparing $\chk{\size}_i(s_j\waf)$ to the computation above. First note that the inversion sequence for $(s_j\waf)^{-1}$ agrees with that of $\waf^{-1}$ with an additional last entry, $\waf(\alpha_j)$. For $j \neq 0$, we have
$    \waf (\alpha_j)
    = w(\alpha_j) + \langle \alpha_j,q\rangle \delta,$
while if $j=0$, then 
$    \waf (-\widetilde{\alpha}+\delta)
    = w(-\widetilde{\alpha}) + (1-\langle \widetilde{\alpha},q\rangle)\delta.$  

When $j\neq i$ this last entry does not affect the computation of $\chk{\size_i}$. Otherwise $\chk{\size}_i(s_j\waf)$ includes the above coefficient of $\delta$ in the sum, but with $\alpha_j$ instead of $\chk{\alpha}_j$. This distinction does not matter if  $\alpha_j$ is a long root since $\alpha_j=\chk{\alpha}_j$ and (similarly for $\widetilde{\alpha}$). But if $\alpha_j$ is short, then difference between the root and coroot in the formulas properly introduces the required scaling factor of $r$.
\end{proof}

\section{Combinatorial Models}
\label{sec:models}


In this section we describe combinatorial models for the affine Weyl groups of classical type, as well as $G_2$, recovering and extending some results in~\cite{hanusa2012abacus,cotton2021core}.  In each case, the models are obtained by exhibiting a suitable equivariant embedding from $\chk{\Q}_{X_n}$ into a type-$A$ coroot lattice $\chk{\Q}_m$. Hence, the objects of these models are partitions, and we interpret the statistic $\chk{\size}$, as well as its refinements $\chk{\size}_i$, in terms of the partitions. Unfortunately, we do not know of such an embedding for the remaining (exceptional) types $F_4$, $E_6, E_7, $ and $E_8$, and so we leave open the problem of finding similar combinatorial models for them.

\begin{figure}[htbp]
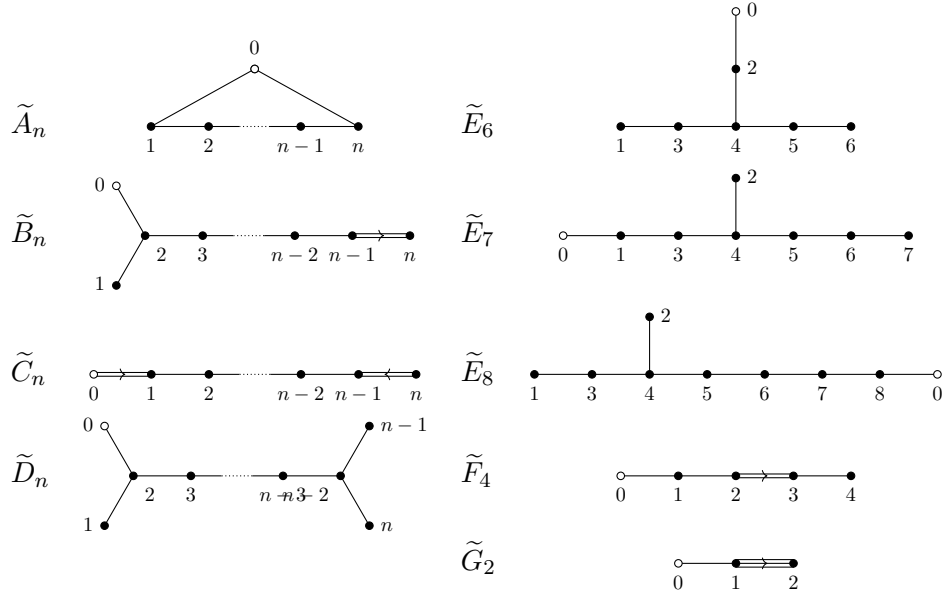

\begin{tabular}{lclc}
$\widetilde{A}_n$ & \dynkin [edge length=2em,extended,labels={0,1,2,n-1,n}]A{}
& $\widetilde{E}_6$ & \dynkin [edge length=2em,extended,labels={0,1,2,3,4,5,6}]E6{}\\
$\widetilde{B}_n$ & \dynkin [edge length=2em,extended,labels={0,1,2,3,n-2,n-1,n}]B{}
& $\widetilde{E}_7$ & \dynkin [edge length=2em,extended,labels={0,1,2,3,4,5,6,7}]E7{}\\
$\widetilde{C}_n$ & \dynkin [edge length=2em,extended,labels={0,1,2,n-2,n-1,n}]C{}
& $\widetilde{E}_8$ & \dynkin [edge length=2em,extended,labels={0,1,2,3,4,5,6,7,8}]E8{}\\
$\widetilde{D}_n$ & \dynkin [edge length=2em,extended,labels={0,1,2,3,n-3,n-2,n-1,n}]D{}
& $\widetilde{F}_4$ & \dynkin [edge length=2em,extended,labels={0,1,2,3,4}]F4{}\\
& &
$\widetilde{G}_2$ & \dynkin [edge length=2em,extended,labels={0,1,2}]G2{}
\end{tabular}
\caption{The affine Dynkin diagrams, with vertex $i$ corresponding to the affine simple reflection $s_i$.}
\label{fig:dynkin}
\end{figure}

\subsection{Type C}
\label{sec:modelC}

The simple roots for $C_n$ are $\alpha_i:=\frac{1}{\sqrt{2}}(e_i-e_{i+1})$ for $1\leq i < n$ and $\alpha_n:=\sqrt 2 e_n$.  Hence, the coroot lattice is $\chk{\Q}_{C_n} = \sqrt{2}\mathbb{Z}^n.$ The action of $\wa(C_n)$ on $\chk{\Q}_{C_n}$ is given explicitly by

\begin{align*}s_i(x_1,\ldots,x_i,x_{i+1},\ldots,x_n) &= (x_1,\ldots,x_{i+1},x_{i},\ldots,x_n) \text{ for } 1 \leq i < n,\\
s_n(x_1,\ldots,x_n) &= (x_1,\ldots,-x_n), \text{ and }\\
s_0(x_1,\ldots,x_n) &= (\sqrt{2}-x_1,\ldots,x_n).
\end{align*}
We embed the type $C_n$ coroot lattice into the coroot lattice for $\widetilde{\mathfrak{S}}_{2n}$ by
\begin{align*}
\iota:\chk{\Q}_{C_n} &\hookrightarrow \chk{\Q}_{2n}, \\
(x_1,\ldots,x_n) &\mapsto \frac{1}{\sqrt{2}} (x_1,\ldots,x_n,-x_n,\ldots,-x_1).
\end{align*}

Evidently, $\iota(\mathbf{x})=\iota(\mathbf{x})^\intercal$ for all $\mathbf{x}\in\chk{\Q}_{C_n}$, and therefore self-conjugate $2n$-cores serve as a combinatorial model for the cores of type $C_n$. It is a particularly well-behaved model because $\iota$ is an isometry (that is, $\langle\iota(\mathbf{x}),\iota(\mathbf{y})\rangle = \langle\mathbf{x},\mathbf{y}\rangle$) and also our definition of $\chk{\size}$ agrees with the number of boxes of the corresponding partitions, as we show in \Cref{thm:cn_cores}. Moreover, it is straightforward to check that the simple reflections of $\wa(C_n)$ agree with the following $\widetilde{\mathfrak{S}}_{2n}$-elements acting on the $\iota$-embedded coroot lattice:
\begin{align*}
s_i &\Leftrightarrow s_i^A s_{2n-i}^A \text{ for } 1 \leq i < n, \text{ while}\\
s_n &\Leftrightarrow s_n^A, \text{ and}\\
s_0 &\Leftrightarrow s_0^A.
\end{align*}

\begin{example}
\label{ex:c2}
The inversion sequence for $\wa({C}_2)$ and $\mathsf{\widetilde{w}}^{-1}=s_1s_0s_2s_1s_0$ (so that $\waf=t_{-\chk{\alpha}_1}s_2$) is
\[-\widetilde{\alpha}+\delta,-\alpha_1-\alpha_2+\delta,-\widetilde{\alpha}+2\delta,-\alpha_2+\delta,-\alpha_1-\alpha_2+2\delta\]
and because $\alpha_1$ is short and $r=2$, we have that $\chk{\size}_1(\mathsf{\widetilde{w}})=2\cdot(1+2)=6$.

On the other hand, observe that $\waf$ corresponds to $q=\sqrt{2}(-1,1) \in \chk{\Q}_{C_2}$.  
 Since $\chk{\omega}_1 = \sqrt{2}(1,0)$ and $c_1=2$, we compute
\[\chk{\size}_1(q) = \left\langle \sqrt{2}(-1,1)-\sqrt{2}(1,0),\sqrt{2}(-1,1)\right\rangle = 2\left\langle (-1,1)-(1,0),(-1,1)\right\rangle = 6.\]

Moreover, the corresponding $4$-core has $6$ boxes with content $1$ or $3$ mod $4$:
$$ \lambda_{(-1,1,1,-1)}=\raisebox{\height}{\young(0123,301,23,1)}. $$
\end{example}

\begin{theorem}
\label{thm:cn_cores}
The map $\mathbf{x}\mapsto \lambda_{\iota(\mathbf{x})}$ is a $\wa(C_n)$-equivarant bijection between the $C_n$ coroot lattice and self-conjugate $(2n)$-cores.  Moreover, for any $\mathbf{x} \in \chk{\Q}_{C_n}$, 
\[
\chk{\size}_i(\mathbf{x})=
\begin{cases}
\chk{\size}_i(\lambda_{\iota(\mathbf{x})}) + \chk{\size}_{2n-i}(\lambda_{\iota(\mathbf{x})}) & \text{if } 1\leq i<n \\ 
\chk{\size}_i(\lambda_{\iota(\mathbf{x})}) & \text{if } i=0,n\end{cases}
,
\] and hence
$\chk{\size}(\mathbf{x})=\chk{\size}(\lambda_{\iota(\mathbf{x})})$.
\end{theorem}

\begin{proof}
Note that the map is well-defined because by definition, $\iota(\mathbf{x})^\intercal=\iota(\mathbf{x})$ and thus by \Cref{prop:ip_content} we have $\lambda_{\iota(\mathbf{x})}$ is self-conjugate. It is evidently injective, and it is easy to see that every self-conjugate partition is in the image as well. Equivariance follows from the straightforward check above and \Cref{thm:sa_eq}. Thus it remains to prove that it preserves the statistic $\chk{\size}$. 

Begin by observing the following:
\[ 
\iota(\chk{\omega_i}) = (\chk{(\omega_i^A)}+\chk{(\omega_{2n-i}^A)}),
\quad 1\leq i<n,
\qquad\text{and}\qquad
\iota(\chk{\omega_i}) = \chk{(\omega_i^A)}, 
\quad i=0,\,n.
\]
where $\chk{(\omega_i^A)}$ is the $i^\text{th}$ fundamental coweight in type $A_{2n-1}$ (and $\chk{(\omega_0^A)}=0$ as always).


Write $\lambda_i$ for the number of boxes in $\lambda_{\iota(\mathbf{x})}$ with content equal to $i$ mod $2n$. The content of \Cref{prop:ip_content} is that $\lambda_i=\left\langle\frac{1}{2}\iota(\mathbf{x})-\chk{(\omega_i^A)}, \iota(\mathbf{x})\right\rangle$. Thus,
for all $1\leq i<n$, we have
\begin{align*}
\chk{\size}_i(\mathbf{x}) &= \left\langle \frac{c_i}{2}\mathbf{x}-\chk{\omega_i}, \mathbf{x} \right\rangle \\
&= \left\langle \mathbf{x}- \chk{\omega_i}, \mathbf{x} \right\rangle \\
&= \left\langle \iota(\mathbf{x}) - \chk{(\omega_i^A)}-\chk{(\omega_{2n-i}^A)}, \iota(\mathbf{x}) \right\rangle \\
&= \left\langle \frac{1}{2} \iota(\mathbf{x}) - \chk{(\omega_i^A)}, \iota(\mathbf{x}) \right\rangle + \left\langle \frac{1}{2} \iota(\mathbf{x}) - \chk{(\omega_{2n-i}^A)}, \iota(\mathbf{x}) \right\rangle \\
&= \lambda_i + \lambda_{2n-i},
\end{align*}
and for $i=0$ or $i=n$ the calculation is similar, but because $c_i=1$, we obtain
\begin{align*}
\chk{\size}_i(\mathbf{x}) 
= \left\langle \frac{1}{2}\mathbf{x}- \chk{\omega_i}, \mathbf{x} \right\rangle 
= \left\langle \frac{1}{2}\iota(\mathbf{x}) - \chk{(\omega_i^A)}, \iota(\mathbf{x}) \right\rangle 
= \lambda_i.
\end{align*}
Finally, summing over all $i$ yields the claim for $\chk{\size}$.

\end{proof}

\subsection{Type B}
\label{sec:modelB}

The simple roots for $B_n$ are $\alpha_i:=e_i-e_{i+1}$ for $1\leq i < n$ and $\alpha_n:=e_n$.  
Hence, the coroot lattice is \[\chk{\Q}_{B_n} = \left\{\mathbf{x} \in \mathbb{Z}^n : \sum_{i=1}^n x_i \equiv 0 \mod 2\right\}.\]  The action of $\wa(B_n)$ on $\chk{\Q}_{B_n}$ is given explicitly by
\begin{align*}s_i(x_1,\ldots,x_i,x_{i+1},\ldots,x_n) &= (x_1,\ldots,x_{i+1},x_{i},\ldots,x_n) \text{ for } 1 \leq i < n,\\
s_n(x_1,\ldots,x_n) &= (x_1,\ldots,-x_n), \text{ and }\\
s_0(x_1,x_2,\ldots,x_n) &= (1-x_2,1-x_1\ldots,x_n).
\end{align*}
We may embed the $B_n$ coroot lattice into the coroot lattice for $\widetilde{\mathfrak{S}}_{2n}$ using essentially the same $\iota$ as for $C_n$, namely $\iota: (x_1,\ldots,x_n) \mapsto (x_1,\ldots,x_n,-x_n,\ldots,-x_1)$. However, because we no longer have the normalization factor, this $\iota$ fails to be an isometry; rather, $\langle\iota(\mathbf{x}),\iota(\mathbf{y})\rangle = 2 \langle\mathbf{x},\mathbf{y}\rangle$.
Nevertheless, it is again straightforward to mimic the action of $\wa(B_n)$ using $\widetilde{\mathfrak{S}}_{2n}$ by
\begin{align*}
s_i &\mapsto s_i^A s_{2n-i}^A \text{ for } 1 \leq i < n, \text{ while}\\
s_n &\mapsto s_n^A, \text{ and}\\
s_0 &\mapsto s_0^As_1^As_{2n-1}^As_0^A.
\end{align*}

We thus obtain a combinatorial model for $\chk{\Q}_{B_n}$.

\begin{theorem}
\label{thm:bn_cores}
The map $\mathbf{x}\mapsto \lambda_{\iota(\mathbf{x})}$ is a $\wa(B_n)$-equivarant bijection between the $B_n$ coroot lattice and self-conjugate $2n$-cores with an even number of boxes on the main diagonal. For $\mathbf{x} \in \chk{\Q}_{B_n}$, 
\[
\chk{\size}_i(\mathbf{x})= \frac{1}{2}
\begin{cases}
\chk{\size}_i(\lambda_{\iota(\mathbf{x})}) + \chk{\size}_{2n-i}(\lambda_{\iota(\mathbf{x})}) & \text{if } 1<i\leq n \\ 
\chk{\size}_1(\lambda_{\iota(\mathbf{x})}) + \chk{\size}_{2n-1}(\lambda_{\iota(\mathbf{x})}) -\chk{\size}_0(\lambda_{\iota(\mathbf{x})}) & \text{if } i=1 \\ \chk{\size}_0(\lambda_{\iota(\mathbf{x})}) & \text{if } i=0
\end{cases}
,
\] and hence
$\chk{\size}(\mathbf{x}) = \frac{1}{2}\left( \chk{\size}(\lambda_{\iota(\mathbf{x})}) - \chk{\size}_0(\lambda_{\iota(\mathbf{x})}) \right)$.
\end{theorem}

\begin{proof}
As in \Cref{thm:cn_cores} we have that $\lambda_{\iota(\mathbf{x})}$ is self-conjugate. Moreover, for any partition $\lambda$, the side length of its Durfee square is the number of black beads that lie below the line in its abacus diagram that witnesses the fact that it is balanced. By definition of $\lambda_{\iota(\mathbf{x})}$, this must be $|x_1|+|x_2|+\cdots+|x_n|$, and since $\mathbf{x}\in \chk{\Q}_{B_n}$, this number must be even. Thus $\lambda_{\iota(\mathbf{x})}$ has an even number of elements along its main diagonal, as desired.

This shows that the map is well-defined. Bijectivity and equivariance follow in a manner analogous to \Cref{thm:cn_cores}. It remains to prove that it has the claimed effect on $\chk{\size}$. We begin as before with the following observation:
\[
\iota(\chk{\omega_i}) = (\chk{(\omega_i^A)}+\chk{(\omega_{2n-i}^A)}),
\quad 1\leq i\leq n,
\qquad\text{and}\qquad
\iota(\chk{\omega_0}) = \chk{(\omega_0^A)}.
\]

Write $\lambda_i$ for the number of boxes in $\lambda_{\iota(\mathbf{x})}$ with content equal to $i$ mod $2n$. By \Cref{prop:ip_content}, this is $\left\langle\frac{1}{2}\iota(\mathbf{x})-\chk{(\omega_i^A)}, \iota(\mathbf{x})\right\rangle$. Thus,
for all $2\leq i<n$,
\begin{align*}
\chk{\size}_i(\mathbf{x}) 
&= \left\langle \mathbf{x}- \chk{\omega_i}, \mathbf{x} \right\rangle \\
&= \frac{1}{2}\left\langle \iota(\mathbf{x}) - \chk{(\omega_i^A)}-\chk{(\omega_{2n-i}^A)}, \iota(\mathbf{x}) \right\rangle \\
&= \frac{1}{2}\left\langle \frac{1}{2}\iota(\mathbf{x}) - \chk{(\omega_i^A)}, \iota(\mathbf{x}) \right\rangle + \frac{1}{2}\left\langle \frac{1}{2}\iota(\mathbf{x}) - \chk{(\omega_{2n-i}^A)}, \iota(\mathbf{x}) \right\rangle \\
&= \frac{\lambda_i + \lambda_{2n-i}}{2}.
\end{align*}

We obtain the calculations for $\chk{\size}_0(\mathbf{x})$ and $\chk{\size}_n(\mathbf{x})$ in a manner analogous to type $C_n$. But when $i=1$, we observe somewhat different behavior because $c_1=1$:
\begin{align*}
\chk{\size}_i(\mathbf{x})
&= \frac{1}{2}\left\langle \frac{1}{2}\iota(\mathbf{x}) - \chk{(\omega_1^A)}-\chk{(\omega_{2n-1}^A)}, \iota(\mathbf{x}) \right\rangle \\
&= \frac{1}{2}\left(\left\langle \frac{1}{2}\iota(\mathbf{x}) - \chk{(\omega_1^A)}, \iota(\mathbf{x}) \right\rangle + \left\langle \frac{1}{2}\iota(\mathbf{x}) - \chk{(\omega_{2n-1}^A)}, \iota(\mathbf{x}) \right\rangle \right. \\
& \qquad\qquad \left. - \left\langle \frac{1}{2}\iota(\mathbf{x}) - \chk{(\omega_0^A)}, \iota(\mathbf{x}) \right\rangle \right) \\
&= \frac{\lambda_1 + \lambda_{2n-1}-\lambda_0}{2}.
\end{align*}
As usual, the last equality follows by applying \Cref{prop:ip_content}.


By comparing the results of these calculations to the definition of $\chk{\size}_i(\lambda_{\iota(\mathbf{x}})$ and summing over $i$, we obtain the desired equalities.
\end{proof}

\subsection{Type D}
\label{sec:modelD}

The simple roots for $D_n$ are $\alpha_i:=e_i-e_{i+1}$ for $1\leq i < n$ and $\alpha_n:=e_{n-1}+e_n$.  The highest root is $\widetilde{\alpha}:=e_1+e_2$, and the coroot lattice is the same as for $B_n$, namely \[\chk{\Q}_{D_n} = \left\{\mathbf{x} \in \mathbb{Z}^n : \sum_{i=1}^n x_i = 0 \mod 2\right\}.\]  The action of $\wa(D_n)$ on $\chk{\Q}_{D_n}$ is given explicitly by
\begin{align*}s_i^D(x_1,\ldots,x_i,x_{i+1},\ldots,x_n) &= (x_1,\ldots,x_{i+1},x_{i},\ldots,x_n) \text{ for } 1 \leq i < n,\\
s_n^D(x_1,\ldots,x_{n-1},x_n) &= (x_1,\ldots,-x_n,-x_{n-1}), \text{ and }\\
s_0^D(x_1,x_2,\ldots,x_n) &= (1-x_2,1-x_1,\ldots,x_n).
\end{align*}
We again embed the $D_n$ coroot lattice into the coroot lattice for $\widetilde{\mathfrak{S}}_{2n}$ using the same $\iota$ as in type $B$ and mimic the action of $\wa(D_n)$ using $\widetilde{\mathfrak{S}}_{2n}$ by
\begin{align*}
s_i &\mapsto s_i^A s_{2n+1-i}^A \text{ for } 1 \leq i < n, \text{ while}\\
s_n &\mapsto s_{n}^A s_{n-1}^A s_{n+1}^A s_n^A, \text{ and}\\
s_0 &\mapsto s_0^As_1^As_{2n-1}^As_0^A.
\end{align*}
We obtain a combinatorial model for $\chk{\Q}_{D_n}$.

\begin{theorem}The map $\mathbf{x}\mapsto\lambda_{\iota(\mathbf{x})}$ is a $\wa(D_n)$-equivarant bijection between the $D_n$ coroot lattice and self-conjugate $2n$-cores with an even number of boxes on the main diagonal. For $\mathbf{x} \in \chk{\Q}_{D_n}$, 
\[
\chk{\size}_i(\mathbf{x})= \frac{1}{2}
\begin{cases}
\chk{\size}_i(\lambda_{\iota(\mathbf{x})}) + \chk{\size}_{2n-i}(\lambda_{\iota(\mathbf{x})}) & \text{if } 1<i<n-1 \\ 
\chk{\size}_{n-1}(\lambda_{\iota(\mathbf{x})}) + \chk{\size}_{n+1}(\lambda_{\iota(\mathbf{x})}) -\chk{\size}_{n}(\lambda_{\iota(\mathbf{x})}) & \text{if } i=n-1 \\
\chk{\size}_1(\lambda_{\iota(\mathbf{x})}) + \chk{\size}_{2n-1}(\lambda_{\iota(\mathbf{x})}) -\chk{\size}_0(\lambda_{\iota(\mathbf{x})}) & \text{if } i=1 \\ \chk{\size}_i(\lambda_{\iota(\mathbf{x})}) & \text{if } i=0,n
\end{cases}
,
\] and hence
$\chk{\size}(\mathbf{x}) = \frac{1}{2}\left( \chk{\size}(\lambda_{\iota(\mathbf{x})}) - \chk{\size}_0(\lambda_{\iota(\mathbf{x})}) - \chk{\size}_n(\lambda_{\iota(\mathbf{x})}) \right)$.
\label{ref:dn_cores}
\end{theorem}

\begin{proof}
Since $\chk{\Q}_B=\chk{\Q}_D$, the first statement is automatic from \Cref{thm:bn_cores}. As usual, to prove the map has the claimed effect on $\chk{\size}$, we begin with the following observation:
\begin{align*}
\iota(\chk{\omega_i}) &= \chk{(\omega_i^A)}+\chk{(\omega_{2n-i}^A)} & 1\leq i<n-1\\
\iota(\chk{(\omega_{n-1}^A)}) &= \chk{(\omega_{n-1}^A)}-\chk{(\omega_n^A)}+\chk{(\omega_{n+1}^A)} \\
\iota(\chk{\omega_n}) &= \chk{(\omega_n^A)}.
\end{align*}

 Therefore, the computation of $\chk{\size}_i$ is identical to the type $B_n$ for all $i$ except for $i=n-1,n$. When $i=n$ it is analogous to the type $C_n$ computation, and when $i=n-1$ we compute:
\begin{align*}
\chk{\size}_{n-1}(\mathbf{x}) 
&= \frac{1}{2}\left\langle \frac{1}{2}\iota(\mathbf{x}) - \chk{(\omega_{n-1}^A)}+\chk{(\omega_{n}^A)}- \chk{(\omega_{n+1}^A)}, \iota(\mathbf{x}) \right\rangle \\
&= \frac{1}{2}\left(\left\langle \frac{1}{2}\iota(\mathbf{x}) - \chk{(\omega_{n-1}^A)}, \iota(\mathbf{x}) \right\rangle - \left\langle \frac{1}{2}\iota(\mathbf{x}) - \chk{(\omega_{n}^A)}, \iota(\mathbf{x}) \right\rangle \right. \\
& \qquad\qquad \left. + \left\langle \frac{1}{2}\iota(\mathbf{x}) - \chk{(\omega_{n+1}^A)}, \iota(\mathbf{x}) \right\rangle \right) \\
&= \frac{\lambda_{n-1} - \lambda_{n} + \lambda_{n+1}}{2}.
\end{align*}
where as usual $\lambda_i$ is the number of boxes in $\lambda_{\iota(\mathbf{x})}$ with content $i$ mod $2n$.
\end{proof}

\subsection{Type G\texorpdfstring{\textsubscript{2}}{2}}
\label{sec:modelG}

Following the usual construction, we consider $G_2$ as acting on the orthogonal complement of the line $\mathrm{span}_\mathbb{R}(1,1,1)$ in $\mathbb{R}^3$.  The simple roots for $G_2$ can be taken to be \[\alpha_1:=(1,-1,0) \text{ and } \alpha_2 := \textstyle\frac{1}{3}(-1,2,-1).\]  With these conventions, 
the coroot lattice is $\chk{\Q}_{G_2} = \chk{\Q}_3$. That is, the type $G_2$ coroot lattice coincides with the coroot lattice for $\widetilde{\mathfrak{S}}_{3}$. Therefore, the map $\mathbf{x}\mapsto\lambda_{\mathbf{x}}$ gives a bijection between $\chk{\Q}_{G_2}$ and $3$-cores.  The action of $\wa(G_2)$ on $\chk{\Q}_{G_2}$ is given explicitly by
\begin{align*}s_1(x_1,x_2,x_3) &= (x_2,x_1,x_3),\\
s_2(x_1,x_2,x_3) &= (-x_3,-x_2,-x_1), \text{ and }\\
s_0(x_1,x_2,x_3) &= (x_3+1,x_2,x_1-1).
\end{align*}
We may therefore emulate the action of $\widetilde{G}_2$ using $\widetilde{\mathfrak{S}}_{3}$ by
\begin{align*}
s_1(\mathbf{x}) &= s_1^A(\mathbf{x}),\\
s_2(\mathbf{x}) &= \mathbf{x}^\intercal, \text{ and}\\
s_0(\mathbf{x}) &= s_0^A(\mathbf{x}).
\end{align*}

As in~\cite{cotton2021core}, we obtain a combinatorial model for $\chk{\Q}_{G_2}$.

\begin{theorem}
The map $q\mapsto \lambda_q$ is a $\wa(G_2)$-equivarant bijection between the $G_2$ coroot lattice and $3$-cores: $s_1^G$ acts on a 3-core by adding or removing all boxes of content $1$, $s_0^G$ acts similarly on boxes of content $0$, and $s_2^G$ acts by conjugation.  For $\mathbf{x} \in \chk{\Q}_{G_2}$,
\[ 
\chk{\size}_i(\mathbf{x}) = 
\begin{cases}
\chk{\size}_0(\lambda_{\mathbf{x}}) & \text{if } i=0 \\ 
\chk{\size}_0(\lambda_{\mathbf{x}}) + \chk{\size}_1(\lambda_{\mathbf{x}}) + \chk{\size}_2(\lambda_{\mathbf{x}}) & \text{if } i=1 \\ 3\cdot \chk{\size}_0(\lambda_{\mathbf{x}}) - \chk{\size}_0(\lambda_{\mathbf{x}}) & \text{if } i=2
\end{cases}
,
\]
and hence 
$\chk{\size}_i(q)=\chk{\size}(\lambda_q)+3\cdot\chk{\size}_2(\lambda_q)$.
\end{theorem}

\begin{proof}
The bijectivity and equivariance of the map is the content of~\Cref{thm:sa_eq}. It remains to prove its effect on $\chk{\size}$.
Because
\[
\chk{\omega_1} = \chk{(\omega_1^A)}+\chk{(\omega_2^A)}
\qquad\text{and}\qquad
\chk{\omega_2} = 3\chk{(\omega_2^A)},
\]
we have
\begin{align*}
\chk{\size}_0(\mathbf{x}) 
&= \left\langle \frac{1}{2}\mathbf{x}- \chk{\omega_0}, \mathbf{x} \right\rangle
= \left\langle \frac{1}{2}\mathbf{x}, \mathbf{x} \right\rangle
&&= \lambda_0
\\
\chk{\size}_1(\mathbf{x}) 
&= \left\langle \frac{3}{2}\mathbf{x}- \chk{\omega_1}, \mathbf{x} \right\rangle
= \left\langle \frac{3}{2}\iota(\mathbf{x}) - \chk{(\omega_1^A)} - \chk{(\omega_2^A)}, \iota(\mathbf{x}) \right\rangle 
&&= \lambda_1+\lambda_2+\lambda_0
\\
\chk{\size}_2(\mathbf{x}) 
&= \left\langle \mathbf{x}- \chk{\omega_2}, \mathbf{x} \right\rangle
= \left\langle \iota(\mathbf{x}) - 3\chk{(\omega_2^A)}, \iota(\mathbf{x}) \right\rangle 
&&= 3\lambda_2-\lambda_0,
\end{align*}
where $\lambda_i$ is as usual the number of boxes in $\lambda_\mathbf{x}$ with content $i$ mod $3$.
\end{proof}

\section{Simultaneous Cores}
\label{sec:sommers}

Recall that an $(a,b)$-core is a partition which is both $a$-core and $b$-core. 
In~\cite[Lemma 3.1]{johnson2015lattice}, Johnson showed that among the set of $a$-cores, the $(a,b)$-cores are precisely those that satisfy certain simple inequalities on the heights of their runners. Thus, when considering them as elements of $\chk{\Q}_{A_{a-1}}$, these inequalities imply that they are lattice points of a polytope in $\mathbb{R}^a$. In fact this polytope is a simplex, previously been considered by Sommers \cite{sommers2005b}, which we describe now. Recall that $\Phi=\Phi(X_n)$ is a root system with irreducible Cartan type $X_n$ and Coxeter number $h$. For $1 \leq i < h$, write $\Phi_{i}$ to denote the set of (positive) roots of height $i$.


\begin{definition}
\label{def:sommers}
For $b$ coprime to $h$, write $\rat=\bt h+\br$ with $\bt,\br\in\Z_{\geq0}$ and $0<\br<h$.  We define the \defn{$b$-Sommers region} 
\begin{align*}
 \Sommers_{X_n}(b):=\left\{x\in V: \begin{array}{rl} \langle x,\alpha\rangle\geq-\bt & \text{for all }\alpha\in \Phi_{\br}\text{ and }\\
 \langle x,\alpha\rangle\leq \bt+1 &\text{for all }\alpha\in \Phi_{h-\br} \end{array} \right\}.
\end{align*}
(We write $\Sommers(b)$ when the root system is clear from context.)
\end{definition}

As in~\Cref{eq:Wcores}, a natural generalization of $\core(a,b)$ to any affine Weyl group is the intersection of the coroot lattice $\chk{\Q}_{X_n}$ with $\Sommers_{X_n}(b)$, so that $\core(a,b) = \core(A_{a-1},b)$.

\subsection{The Sommers Region and the Fundamental Alcove}

We would like to perform the $\chk{\size}$-weighted enumeration of $\core(X_n,b)$ using Ehrhart theory. Unfortunately, the family \[\{\Sommers_{X_n}(b): \gcd(b,h)=1\}\] does not consist of dilations of a fixed polytope---but this difficulty can be circumvented.
Define, for any $x\in V$, the statistic
\begin{equation}\sizeb(x):=\frac{h}{2}\left(\left\|x - \frac{b \chk{\rho}}{h}\right\|^2 - \left\|\frac{\chk{\rho}}{h}\right\|^2 \right). \end{equation}
Notice that when $q\in\chk{\Q}$, $w\in W$, and $b=1$, we have $\size^{(1)}(w(q))=\chk{\size}(w(q)) = \chk{\size}(wt_q)$ (recalling \Cref{thm:sizer}).

We recall from~\cite[\S 4]{thiel2017strange} that there is a unique element $\wf\in \wa$ such that $\frac{b}{h}\chk{\rho} = \wf(\frac{\chk{\rho}}{h})$, and that left-multiplication by this element maps $\Sommers(b)$ onto the $b$-fold dilation of the fundamental alcove $\ac:= \{x\in V:\langle x,\alpha\rangle\geq 0\text{ for all }\alpha\in\simp\text{ and }\langle x,\tilde{\alpha}\rangle\leq 1\}$. It also respects the lattice points in the following sense:

\begin{theorem}
\label{thm:size_after_wb}
For $b$ coprime to $h$, the following holds as an equality of multisets:
\[\Big\{\chk{\size}(q) : q \in \core(X_n,b)\Big\} = \Big\{\sizeb(q) : q \in \rat\ac \cap \chk{\Q}_{X_n}\Big\}.\]
\end{theorem}
\begin{proof}
We first note that since $\wf$ maps $\Sommers(b)$ onto $b\ac$, and also $\wf\in\wa$ and thus is a $\chk{\Q}$-preserving bijection, it restricts to a bijection $\core(\wa,b)\to \rat\ac \cap \chk{\Q}$.  Write $\wf = w t_{q_0}$; then $\wf = t_{\frac{b\chk{\rho}}{h}} w t_{-\frac{\chk{\rho}}{h}}$. Since $\|\cdot\|$ is $W$-invariant, for $q\in \core(\wa,b)$:
\begin{align*}
\sizeb(\wf(q)) &= \frac{h}{2}\left(\left\| t_{\frac{b}{h}\chk{\rho}} w t_{\frac{1}{h}\chk{\rho}} (q) - \frac{b \chk{\rho}}{h}\right\|^2 - \left\|\frac{\chk{\rho}}{h}\right\|^2 \right) \\
&= \frac{h}{2}\left(\left\|w\left(q - \frac{\chk{\rho}}{h}\right)\right\|^2 - \left\|\frac{\chk{\rho}}{h}\right\|^2 \right) = \chk{\size}(q).\qedhere
\end{align*}
\end{proof}

\subsection{Maximum Size}

\Cref{thm:size_after_wb} is a primary tool in computing the expected size of simultaneous cores in the next section. As a simpler application, we proceed as in \cite{thiel2017strange} to determine the maximum $\chk{\size}$ of a simultaneous core (extending that result to the non-simply-laced types).

\begin{theorem}
\label{thm:max_size}
For $\widetilde{W}$ an irreducible affine Weyl group with $\gcd(h,b)=1$, 
\[
\max{q \in \Sommers(b)\cap \chk{\Q}}{\chk{\size}(q)} = \frac{r\chk{g}}{h}\frac{n(b^2-1)(h+1)}{24}.
\] 
Moreover, this maximum is attained by a unique point $q_* \in \Sommers(b)$.
\end{theorem}

\begin{proof}
We claim that the maximum is obtained at $q_* = \wf^{-1}(0)$, where $\wf\in\widetilde{W}$ is the same element as used in \Cref{thm:size_after_wb}. First note that since $\wf$ maps $\Sommers(\rat)$ bijectively to $b\ac$, this $q_*$ is indeed in $\core(\wa,\rat)=\Sommers(\rat)\cap\Q$.
Since $\wf$ maps $\chk{\size}$ to $\size^{(b)}$, we will show the equivalent statement that $0$ is the unique element of $\rat\ac\cap\chk\Q$ of maximum $\size^{(b)}$.

Since $\left\|\frac{\chk{\rho}}{h}\right\|^2$ is a constant, it suffices to maximize $Q(\mathbf{x})=\left\|x-\rat\rhoh\right\|^2$. But the fact that $0$ maximizes $Q(\mathbf{x})$ over $\rat\ac$ is known; it follows for instance from the ``very strange'' formula of Kac (cf.~\cite[Equation (0.9)]{kac84theta}). Moreover, $Q(\mathbf{x})$ is a strictly convex function, so it can only be maximized at a vertex of the convex polytope $\rat\ac$. However, no other vertices of $\rat\ac$ are in $\chk\Q$, which implies that $0$ is the only point in $\rat\ac\cap\chk\Q$ of maximum $\size^{(b)}$.\\

Moreover, we may explicitly compute $\size(q_*)$ as follows.
\[
\size^{(b)}(0) 
=\frac{h}{2}\left\|\rat\rhoh\right\|^2-\frac{h}{2}\left\|\rhoh\right\|^2
=(b^2-1)\frac{1}{2h}\left\langle\chk{\rho},\chk{\rho} \right\rangle.
\]

The desired formula then follows from the explicit computation of $\langle \chk{\rho},\chk{\rho}\rangle$, a dual version of the ``strange formula'' of Freudenthal and de Vries:

\begin{theorem}[see {\cite[Section 4]{burns2000elementary}}]
\label{thm:strange_dual}
For $\widetilde{W}$ an irreducible affine Weyl group,
\[\langle \chk{\rho},\chk{\rho}\rangle = r\chk{g}\frac{n(h+1)}{12}.\qedhere\]
\end{theorem}
\end{proof}

Although~\Cref{thm:max_size} proves that $\chk\size(\wf^{-1}(0))$ is the maximum that the statistic $\chk\size$ can take on $\core(\wa,b)$, more is true in type $A_n$, where J.~Vandehey shows that the largest $(a,b)$-core contains all other $(a,b)$-cores as subdiagrams (see~\cite{vandehey2008general,fayers2011t}). Previously the second and third authors conjectured \cite[Conjecture 6.14]{thiel2017strange} that the inversion set of $\wf$ contains the inversion sets of all other affine elements corresponding to elements of $\core(\wa,b)$. We here extend that conjecture to the non-simply-laced types as well.

\begin{conjecture}
The element $\wf$ is maximal in the weak order on $\wa/W$ among all dominant elements $\{ \waf \in \wa/W : \waf^{-1}(0) \in \Sommers(b)\}$.
\label{conj:w_b_is_maximal_in_weak}
\end{conjecture}

\subsection{Simultaneous Combinatorial Models}

One of the major advantages to our improvement from simply-laced to general affine Weyl groups is the ability to incorporate type $C_n$, whose cores lying in $\Sommers_{C_n}(b)$ also have a natural combinatorial model:

\begin{theorem}
The map $\mathbf{x}\mapsto \lambda_{\iota(\mathbf{x})}$ of \Cref{thm:cn_cores} restricts to a bijection from $\Sommers_{C_n}(b)\cap\chk{\Q}_{C_n}$ to self-conjugate $(2n,b)$-cores.
\end{theorem}

\begin{proof}
Recall that a $2n$-core is self-conjugate if and only if $q_\lambda = (q_\lambda)^\intercal$. In type $C_n$, the roots of height $r=2k-1$ are $\frac{1}{\sqrt{2}}(e_i-e_{i+r})$ for all $i\leq n-r$, as well as $\frac{1}{\sqrt{2}}(e_{n-i+1} + e_{n-r+i})$ for $1\leq i\leq k$. Write $\iota(\mathbf{x})=(y_1,\dots y_{2n})$.

Thus  $\langle\mathbf{x},\alpha\rangle \geq -t $ for all $\alpha\in\Phi_r(C_n)$ if and only if 
\begin{align*}
y_i-y_{i+r} &\geq -t & 1\leq i\leq n-r \\
-y_{2n+1-i}+y_{2n+1-i-r} &\geq -t & 1\leq i\leq n-r \\
y_{n-i+1}+y_{n-r+i} &\geq -t & 1\leq i\leq \textstyle\frac{r+1}{2} \\
-y_{n+i}-y_{n+1+r-i} &\geq -t & 1\leq i\leq \textstyle\frac{r+1}{2}
\end{align*}
By definition of $\iota_C$ we have that $y_j=-y_{2n+1-j}$ for all $j$. Thus under the assumption that $\mathbf{y}=\mathbf{y}^\intercal$, the above inequalities are precisely the system $y_i-y_{i+r}\geq -t$ for all $1\leq i\leq 2n-r$, which is to say, $\langle\iota(\mathbf{x}),\alpha\rangle \geq -t$ for all $\alpha\in \Phi_r(A_{2n-1})$. In a similar way, under the same assumption $\langle\mathbf{x},\alpha\rangle \leq t+1$ for all $\alpha\in \Phi_{n-r}(C_n)$ is equivalent to $\langle\iota_C(\mathbf{x}),\alpha\rangle \leq t+1$ for all $\alpha\in \Phi_{n-r}(A_{2n-1})$.

Therefore, the target of our map, restricted to $\Sommers_{C_n}(b)\cap \chk{\Q}_{C_n}$, is indeed $\Sommers_{A_n}(b)\cap \chk{\Q}_{A_{2n-1}}$. Moreover, the map $(y_1,\dots, y_{2n})\mapsto \frac{1}{\sqrt{2}}(y_1,\dots, y_n)$ is a well-defined inverse.
\end{proof}

The above theorem permits us to understand the left-hand side of \Cref{thm:max_size} and \Cref{thm:main_thm} as the maximum and expected $\chk{\size}$ of self-conjugate simultaneous cores.  In particular, this thus recovers ``half'' of the Chen--Huang--Wang result (that is, the case when $a$ is even). 

Unfortunately, this combinatorial interpretation appears to be limited to type $C_n$.  Even for the other classical types, the maps $\iota$ from \Cref{sec:size_gen} do not map the Sommers regions into type-$A$ Sommers regions. Since there are many possible embeddings of $\chk{\Q}_{X_n}$ into $\chk{\Q}_m$ for various $m$, it is possible that such a deficiency may be overcome. It would already be interesting to understand combinatorial conditions on the $3$-cores that do lie in the image of the $G_2$ Sommers region.

\section{Expected Size of Simultaneous Cores}
\label{sec:expectation}

We are now ready to prove \Cref{thm:main_thm}:

\mainthm*

We do this by computing the left-hand side explicitly for each type $X_n$, which by \Cref{thm:size_after_wb} is
\begin{align*}
    \Expt{q \in \core(X_n,b)}{\chk{\size}(q)} 
    &= \frac{1}{|\core(X_n,b)|}\sum_{q\in \core(X_n,b)} \chk{\size}(q) = \frac{1}{|b\ac\cap\chk{\Q}|}\sum_{q\in b\ac\cap\chk{\Q}} \sizeb(q).
\end{align*} 
The denominator was explicitly and nearly-uniformly calculated by Haiman \cite{haiman1994conjectures}. To compute the sum, we first record the vertices of the fundamental alcove $\ac$: they are \(\Gamma:=\{0\}\cup\left\{\frac{\chk{\omega}_i}{c_i} : 1 \leq i \leq n\right\}\), where we recall that the $c_i$ are defined by $\widetilde{\alpha}=\sum_{i=1}^nc_i\alpha_i$.  As in \cite{thiel2017strange}, we proceed by translating the problem to the coweight lattice.  Define the \defn{extended affine Weyl group} by $\wex:=W\ltimes \cwl$, and write the group of automorphisms for $b\ac$ as $b\cycl:=\{\waf\in\wex:\waf(b\ac)=b\ac\}$. These groups are isomorphic for all $b$, and in particular have constant order that we denote by $f$.

\begin{proposition}[{\cite[Theorem 2.5 \& Lemma 6.11]{thiel2017strange}}]
For any $b$ coprime to the Coxeter number $h$ in type $X_n$:
	\begin{enumerate}[(a)]
	\item The action of $b\cycl$ on $b\ac\cap\cwl$ is free.
	\item Each $b\cycl$ orbit of $b\ac\cap\cwl$ contains exactly one element of $b\ac\cap\chk{\Q}=\core(X_n,b)$.
	\item For any $\waf\in b\cycl$ and any $\omega\in\cwl$, $\size^{(b)}(\omega)=\size^{(b)}(\waf\cdot\omega)$
	\end{enumerate}
\label{prop:barational}
\end{proposition}

Using this, we finish the translation from $\chk{\Q}$ to $\cwl$:
\begin{equation}
\label{eq:ev_coweight_sum}
\Expt{q \in \core(X_n,b)}{\chk{\size}(q)} = \frac{1}{|b\ac\cap\chk{\Q}|}\cdot \frac{1}{f}\sum_{q\in b\ac\cap\cwl} \sizeb(q).    
\end{equation}

We now recall the relevant Ehrhart-theoretic tools. For any degree-$r$ polynomial $F:\R^n\to\R$, its \defn{weighted lattice point enumerator} over $\cwl$ is
$ \ac_F(b) := \sum\limits_{q\in b\ac\cap\cwl} F(q)$. This $\ac_F(b)$ is a quasipolynomial in $b$, of degree $n+r$ and period $c=c(\widetilde{X}_n):=\lcm(c_1,\dots,c_n)$, where the $c_i$ are again the denominators of the vertices of $\ac$. As $\size^{(b)}$ changes with $b$, Ehrhart theory appears to be inapplicable---however, a judicious rewriting shows that this is not the case.

\begin{proposition}
	The weighted lattice point enumerator $\ac_{\sizeb}(b)$
    is a quasipolynomial in $b$ of degree $n+2$ and period       $c(\widetilde{X}_n):=\lcm(c_1,\ldots,c_n)$.
\label{prop:quasipoly}
\end{proposition}

\begin{proof}
    Notice that $\chk{\size}_b(x)=\frac{h}{2}\left\|x\right\|^2-b\langle x,\chk{\rho}\rangle+(b^2-1)\frac{\left\|\check{\rho}\right\|^2}{2h}.$
   Thus we find that
   \[\ac_{\sizeb}(b)={\textstyle\frac{h}{2}}\ac_{\|\cdot\|^2}(b)-b\ac_{\langle\cdot,\chk{\rho}\rangle}(b)+(b^2-1)\ac_{\frac{\left\|\chk{\rho}\right\|^2}{2h}}(b)\]
   is a quasipolynomial in $b$ of degree $n+2$ and period $c(\widetilde{X}_n)$.
\end{proof}

Therefore, to complete the proof of \Cref{thm:main_thm}, we may compute the quasipolynomial on for all components that contain a residue $b\mod c_{X_n}$ that is coprime to $h$. For the exceptional types, this is already a finite and computationally feasible calculation. 

Types $A$ and $D$ are simply-laced, and so their calculation was already completed in \cite{thiel2017strange}. Thus, we proceed along similar lines for types $B$ and $C$. In what follows, we write $\ac_{\sizeb}(b)_i$ to mean the polynomial which agrees with $\ac_{\sizeb}(b)$ for all $b\equiv i\bmod c(\widetilde{X}_n)$.  Relevant data to complete these computations for the irreducible root systems is provided in~\Cref{fig:data}.

\begin{figure}[htbp]
\[\begin{array}{c|c|c|c|c|c|c|c}
X_n & h & \chk{g} & e_i & c_i & f & r & m(X_n) \\ \hline
A_n & n+1 & n+1 & 1,2,\ldots,n & 1,1,\ldots,1,1 & n+1 & 1 & 1 \\
B_n & 2n & n+1 & 1,3,\ldots,2n-1 & 1,2,\ldots,2,2 & 2 & 2 & 2\\
C_n & 2n & 2n-1 & 1,3,\ldots,2n-1 & 2,2,\ldots,2,1 & 2 & 2 & 2\\
D_n & 2n-2 & 2n-2 & 1,3,\ldots,2n-3,n-1 & 1,2,\ldots,2,1,1 & 4 & 1 & 2 \\
E_6 & 12 & 12 & 1, 4, 5, 7, 8, 11 & 1, 2, 2, 3, 2, 1 & 3 & 1 & 6 \\
E_7 & 18 & 18 & 1, 5, 7, 9, 11, 13, 17 & 2, 2, 3, 4, 3, 2, 1 & 2 & 1 & 12 \\
E_8 & 30 & 30 & 1, 7, 11, 13, 17, 19, 23, 29 & 2, 3, 4, 6, 5, 4, 3, 2 & 1 & 1 & 60 \\
F_4 & 12 & 9 & 1,5,7,11 & 2, 3, 4, 2 & 1 & 2 & 12 \\
G_2 & 6 & 4 & 1,5 & 3,2 & 1 & 3 & 6 \\
\end{array}\]
\caption{The type $X_n$, Coxeter number $h$, dual Coxeter number $\protect\chk{g}$, exponents $e_i$, coefficients of the highest root $c_i$, index of connection $f$, ratio of long to short root $r$, and $\mathrm{gcd}$ of the $c_i$ for the irreducible root systems.}
\label{fig:data}
\end{figure}

\subsection{Types B and C}

For either $X_n=B_n$ or $X_n=C_n$ we have $c(\widetilde{X}_n)=2$ and exponents $1,3,\ldots,2n-1$. Since all $b$ in the theorem statement are coprime to $h$, they must be odd, and so it suffices to compute only the polynomial $\ac_{\sizeb}(b)_1$. This polynomial has degree $n+2$. Recall the fact that the interior of $e_i\ac$ contains no points of $\chk{\Lambda}_{X_n}$ for all exponents $e_i$, because they are strictly less than $h$ (see e.g. \cite[Section 7.4]{thiel2017strange}). 
Thus, by Ehrhart reciprocity we conclude that $\ac_{\sizeb}(-e_i)=0$, and because in types $B_n$ and $C_n$ the exponents are the odd integers $1,3,\dots, 2n-1$, we thus we have $n$ roots of our desired polynomial. 

Moreover, it is easy to compute that $\ac_{\size}(1)=0$, and therefore (c.f. \cite[Proposition 7.5]{thiel2017strange}) we also have $\ac_{\sizeb}(-h-1)=0$. Hence we know all $n+2$ of our desired polynomial's roots, that is:
$$\ac_{\sizeb}(b)_1 = \kappa(b-1)(b+2n+1)\prod_{j=1}^{n} (b+e_j).$$
where the leading coefficient $\kappa$ depends only on $n$ and whether we are in types $B_n$ or $C_n$. Importantly, this formula holds for any odd $b$, not merely those $b$ which are coprime to $h=2n$. Thus we can determine it by explicitly calculating $\ac_{\sizeb}(b)$ at $b=3$, as follows:
\[ 
\kappa = \frac{\ac_{\size^{(3)}}(3)}{(3-1)(3+2n+1)\prod_{j=1}^{n} (3+(2j-1))} = \frac{1}{2^{n+2}(n+2)!}\ac_{\size^{(3)}}(3).
\]

From here, we recall Haiman's enumeration \cite{haiman1994conjectures}:  $|b\ac\cap\chk{\Q}|=\frac{1}{|W|}\prod_{j=1}^n(b+e_j)$, whenever $b$ and $h$ are coprime. Thus, we rewrite \Cref{eq:ev_coweight_sum} as the quadratic polynomial:
\begin{align*}
\Expt{q \in \core(X_n,b)}{\chk{\size}(q)} &= \frac{2^n n!\kappa}{f}(b-1)(b+2n+1) \\
&= \frac{2^n n!}{2\cdot 2^{n+2}(n+2)!}\ac_{\size^{(3)}}(3)\cdot (b-1)(b+h-1) \\
& = \ac_{\size^{(3)}}(3)\cdot \frac{(b-1)(b+h-1)}{8(n+1)(n+2)}
\end{align*}
Thus, to prove the formula from \Cref{thm:main_thm}, what remains to be verified is that
\[
\ac_{\size^{(3)}}(3)=\frac{r\chk{g}}{3h}\cdot n(n+1)(n+2) =
\begin{cases}
\frac{1}{3}(n+1)^2(n+2) & \text{in type } B_n, \\
\frac{1}{3}(n+1)(n+2)(2n-1) & \text{in type } C_n.
\end{cases}
\]

\subsubsection{Type B}
\label{sec:type_b_calc}

Beginning in type $B_n$, observe that since $\sum_{j=1}^n a_j \chk{\omega}_j = (\sum_{j=1}^n a_j,\sum_{i=2}^n a_j,\ldots,a_n)$ and $\widetilde{\alpha}=(1,1,0,0,\ldots,0)$, the coweight points in $3\ac$ are those points for which $a_1+2(\sum_{j=2}^n a_j) \leq 3$, which are the $2n+2$ points \[\{0,\chk{\omega}_1,2\chk{\omega}_1,3\chk{\omega}_1\}\cup \{\chk{\omega}_j,\chk{\omega}_1+\chk{\omega}_j\}_{j=2}^n.\]
We compute that for $2 \leq j \leq n$
\begin{align*}\size^{(3)}(\chk{\omega}_j) &= \frac{n}{(2n)^2}\left( \left(\|2n\chk{\omega}_j-3\chk{\rho}\|\right)^2-\|\chk{\rho}\|^2\right)
\\ &=
\frac{1}{4n} \left( \sum_{i=1}^j \big(2n-3(n+1-i)\big)^2+\sum_{i=j+1}^n \big(3(n+1-i)\big)^2 -\frac{n(n+1)(2n+1)}{6}\right) \\ 
 &= \frac{(2n+2-3 j) (2n+1-3 j)}{6}. 
\end{align*}
A similar computation also shows that $\size^{(3)}(\chk{\omega}_1+\chk{\omega}_j) = \frac{(2n+2-3 j) (2n+1-3 j)}{6}$
and therefore 
\[\sum_{j=2}^n \size^{(3)}(\chk{\omega}_j)+\size^{(3)}(\chk{\omega}_1+\chk{\omega}_j) = \frac{(n-2) (n-1)^2}{3}.\]
The remaining points give contributions of the following form: 
\begin{align*}
\size^{(3)}(k\chk{\omega}_1) &= \frac{n}{(2 n)^2} \left( (2 k n - 3n)^2 + 
     \sum_{i=2}^n (3 (n + 1 - i))^2 - \sum_{i=1}^n (n + 1 - i)^2 \right) \\
    &=\frac{1 + 3 n - 9 k n + 3 k^2 n + 2 n^2}{3}.
\end{align*}

We thus compute that $\sum_{k=0}^3 \size^{(3)}(k\chk{\omega}_1)= \frac{4(1 + 2 n^2)}{3}$. Putting these together, we obtain 
$$ \ac_{\size^{(3)}}(3) = \frac{{n-2}{n-1}^2}{3} + \frac{4(1+2n^2)}{3} = \frac{1}{3}(n+1)^2(n+2) $$
which concludes the proof of \Cref{thm:main_thm} for type $B_n$.

\subsubsection{Type C}
\label{sec:type_c_calc}
Turning to type $C_n$, we have $\widetilde{\alpha}=\frac{2}{\sqrt{2}}e_1$ and
$$\sum_{j=1}^n a_j \chk{\omega}_j = \sqrt{2}\left(\frac{a_n}{2}+\sum_{j=1}^{n-1} a_j,~\frac{a_n}{2}+\sum_{j=2}^n a_j,\ldots,~\frac{a_n}{2}+a_{n-1},~\frac{a_n}{2}\right)$$
So the coweight points in $3\ac$ are those points for which $2(\sum_{j=1}^{n-1} a_j)+a_n \leq 3$, which are the $2n+2$ points

\[\{0,\chk{\omega}_n,2\chk{\omega}_n,3\chk{\omega}_n\}\cup \{\chk{\omega}_j,\chk{\omega}_n+\chk{\omega}_j\}_{j=1}^{n-1}.\]
Using the fact that $\chk{\rho}=\sum_{j=1}^n \chk{\omega}_j = \frac{\sqrt{2}}{2}(2n-1,2n-3,2n-5,\ldots,1)$, we compute that for $1 \leq j \leq n-1$
\begin{align*}
\size^{(3)}(\chk{\omega}_j) &= \frac{n}{(2n)^2}\left( \left(\|2n\chk{\omega}_j-3\chk{\rho}\|\right)^2-\|\chk{\rho}\|^2\right) \\ &=
\frac{1}{4n}\left( \sum_{i=1}^j 2\big(2n-\frac{3}{2}(2n+1-2i)\big)^2+\sum_{i=j+1}^n \frac{(3(2n+1-2i)\big)^2}{2} -\frac{4n^3-n}{6}\right) \\ 
 &= \frac{(2n-3j)^2-1}{3}.
\end{align*}
A similar computation also shows that $\size^{(3)}(\chk{\omega}_n+\chk{\omega}_j) = \frac{(n-3j)^2-1}{3}$ and therefore
\[\sum_{j=1}^{n-1} \size^{(3)}(\chk{\omega}_j)+\size^{(3)}(\chk{\omega}_n+\chk{\omega}_j) = \frac{(n-2) (n-1)(2n+1)}{3}.\]
The remaining points give contributions of the following form:
\begin{align*}
\size^{(3)}(k\chk{\omega}_n) &= \frac{n}{(2 n)^2} \left(
     \sum_{i=1}^n 2 \left(kn-\frac{3}{2}(2n+1-2i)\right)^2 -\sum_{i=1}^n \frac{(2n+1-i)^2}{2} \right) \\
    &=\frac{-2+8 n^2-9 k n^2+3 k^2 n^2}{6}.
\end{align*}

We thus compute that $\sum_{k=0}^3 \size^{(3)}(k\chk{\omega}_1)  = \frac{2(5n^2-2)}{3}$. Putting these together, we obtain
$$ \ac_{\sizeb}(b) = \frac{1}{3}(n+1)(n+2)(2n-1),$$
which concludes the proof of \Cref{thm:main_thm} for type $C_n$.

\subsection{Types F\texorpdfstring{\textsubscript{4}}{4} and G\texorpdfstring{\textsubscript{2}}{2}}
\label{sec:type_fg_calc}

Finally, we return to the exceptional types. Since we had previously handled the simply-laced types $E_6, E_7,$ and $E_8$ \cite[Section 7.6]{thiel2017strange}, it remains only to confirm our formula for $F_4$ and $G_2$. Since $\ac_{\sizeb}(b)$ is a polynomial in each residue class mod $m(\widetilde{X}_n)$ (here, 12 for $F_4$ and 6 for $G_2$) of degree $n+2$, we can simply compute $\ac_{\sizeb}(b)$ for enough values of $b$ in each relevant residue class and perform Lagrange interpolation. For all relevant residue classes in both groups there is exactly one exponent $e_j$ and so ``enough'' values of $b$ means $6$ for $F_4$, and $4$ for $G_2$.

When performing this computation in \texttt{SAGE}, we see that the polynomials $\ac_{\sizeb}(b)_i$ coincide for the relevant residue classes. Namely:
\[
\ac_{\sizeb}(b)=
\begin{cases}
\frac{1}{18432}(b-1)(b+1)(b+5)(b+7)(b+11)(b+13) & \text{in type } F_4, \\
\frac{1}{144}(b-1)(b+1)(b+5)(b+7) & \text{in type } G_2.
\end{cases}
\]
for all $b$ coprime to $h$ (that is, in both cases, for $b\equiv 1,5\bmod 6$). Dividing these polynomials by $f$ and $|b\ac\cap\chk{\bt}|$ as in \Cref{eq:ev_coweight_sum}, we obtain a formula for the expected $\chk{\size}$ of a core that agrees with
\[
\Expt{q \in \core(X_n,b)}{\chk{\size}(\lambda_q)}=\frac{r \chk{g}}{h}\frac{n(b-1)(h+b+1)}{24}.
\]
This completes the proof of \Cref{thm:main_thm} for $F_4$ and $G_2$, and thus for all types.

\section*{Acknowledgements} 
We thank Benjamin Cotton for help drawing~\Cref{fig:a2}. The first author was partially supported by NSF grant 1601961, NSF grant 1745638, and Czech Science Foundation grant 21-00420M. The third author was partially supported by NSF grant 2246877.


\bibliography{strange2full}{}
\bibliographystyle{amsalpha}

\end{document}